\documentclass[peerreview]{article}
\usepackage[utf8]{inputenc}
\usepackage{amsmath}
\usepackage{amssymb}
\usepackage{commath}
\usepackage{multirow}
\usepackage{algorithm}
\usepackage{algpseudocode}
\usepackage{subcaption}
\usepackage{float}
\usepackage{graphicx}
\usepackage[caption = false]{subfig}
\usepackage[english]{babel}
\usepackage[utf8]{inputenc}
\usepackage{adjustbox}
\usepackage{url}
\usepackage{outlines}
\usepackage{authblk}
\usepackage{titling}
\usepackage{ltablex}      % For longtable + tabularx functionality
\usepackage{booktabs}     % For better table lines (optional but recommended)
\keepXColumns              % Keeps X columns consistent across pages

\usepackage{longtable}
\usepackage{rotating}
\usepackage{adjustbox}
\usepackage{geometry}
\usepackage{placeins}
\usepackage{siunitx}           % optional: for units
\usepackage{caption}
\usepackage{makecell}

\keepXColumns % Optional: keeps column widths consistent across page breaks

\usepackage[nottoc]{tocbibind}

\usepackage[style=ieee, backend=biber]{biblatex} 
\title{Ensemble Parameter Estimation for the Lumped Parameter Linear Superposition (LPLSP) Framework: A Rapid Approach to Reduced-Order Modeling for Transient Thermal Systems}
\author{Neelakantan Padmanabhan}
\affil{ZF Active Safety Systems US Inc, \\ Livonia, USA, \\ \textit{neel.padmanabhan@zf.com}}
\author{}
\affil{}
\date{\today}

\begin{document}
\maketitle
\section{Abstract}
Efficient thermal modeling is essential for the design and reliability of power electronics systems, particularly under fast transient operating conditions. Building upon our previous Lumped Parameter Linear Superposition (LPLSP) approach, this work introduces an ensemble parameter estimation framework that enables reduced order thermal model generation from a single transient dataset. Unlike the earlier implementation that relied on multiple parametric simulations to excite each heat source independently, the proposed approach simultaneously identifies all model coefficients using fully transient excitations. Two estimation strategies namely rank-reduction and two-stage decomposition are developed to further reduce computational cost and improve scalability for larger systems. The proposed strategies yield ROMs with mean temperature-prediction errors within $5\%$ of CFD simulations while reducing model-development times to $O(10^0 s)$--$O(10^1 s)$. Once constructed, the ROM evaluates new transient operating conditions in $O(10^0 s)$, enabling rapid thermal analysis and enabling automated generation of digital twins for both simulated and physical systems. These advancements significantly accelerate early-stage design iterations and mission-profile evaluations for thermal management of electronic systems.
%-----------------
\section{Introduction}
\begin{figure}[!t]
	\centering
	\includegraphics[width=0.9\linewidth]{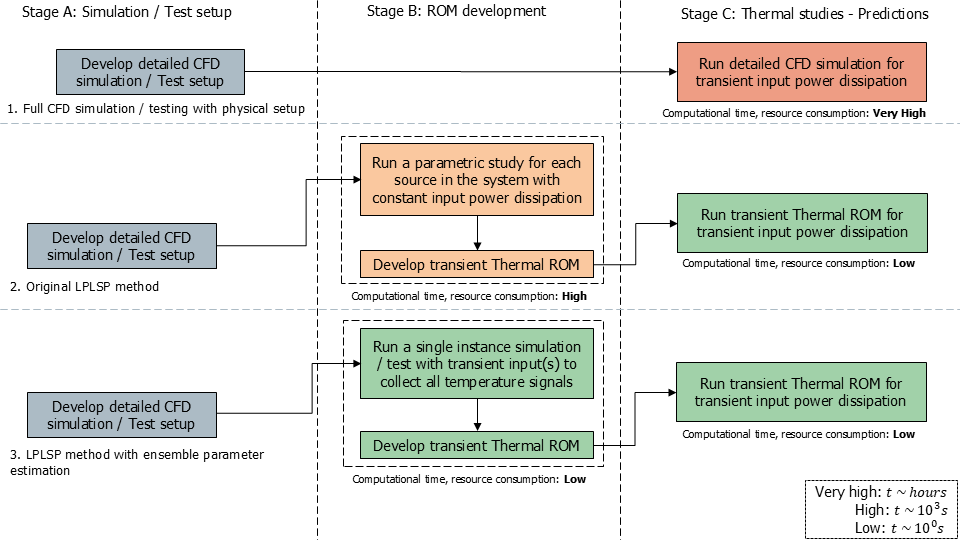}
	\caption{Stages in the development of a reduced-order / computational model. The reported time scales correspond to execution times on a workstation equipped with an Intel\textsuperscript{\textregistered} Core\textsuperscript{\texttrademark} i7--10850H processor (12~cores, 2.7~GHz) and 32~GB RAM.}
	\label{fig:ROMstages}
\end{figure}

Reduced-order modeling (ROM) has become an essential strategy for accelerating transient simulations in thermal management of electronic systems, particularly when full high-fidelity CFD analyses \cite{fullSim1, fullSim2, fullSim3, neel_ecce, neel_jot} are prohibitively expensive for iterative design or mission-profile evaluations. Existing approaches for ROM development include Compact Thermal Networks (CTNs) \cite{CTN1,CTN2,CTN3,CTN4,CTN5}, projection-based techniques such as Krylov subspace and Proper Orthogonal Decomposition (POD) \cite{krylov1,krylov2,krylov3,krylov4,krylov5,BT1,BT2, neel_chapter,POD1,POD2,POD3,POD4,POD5,POD6}, and data-driven methods including physics-informed neural networks \cite{ML1,ML2,ML3,PINN1,PINN2,PINN3,PINN4, OPIN1, OPIN2}. These methods offer varying trade-offs between accuracy, interpretability, and computational cost. CTNs provide physically interpretable models but often require extensive measurements and struggle with multi-source coupling. Projection-based ROMs deliver strong performance in linear regimes but scale poorly for nonlinear, time-varying systems. Machine learning approaches promise flexibility and have gained traction for solving high-dimensional PDEs by embedding governing equations into neural architectures, enabling surrogate models that generalize across parametric variations. However, they often require large, high-quality datasets and extensive hyperparameter tuning, which can be prohibitive from computation time and cost perspective for industrial workflows. Moreover, PINNs and purely data-driven surrogates lack formal guarantees of stability and robustness under extrapolation, raising concerns for safety-critical applications such as automotive power electronics. A detailed review of these methods and their limitations was presented in our previous work \cite{LPLSP2}, and therefore is not repeated here.
Building on this foundation, the Lumped Parameter Linear Superposition (LPLSP) method introduced in \cite{LPLSP1,LPLSP2} demonstrated that accurate transient temperature prediction can be achieved at a fraction of the computational cost of CFD, making it highly attractive for rapid thermal evaluation and design exploration. While the method enables rapid temperature prediction for new input conditions, the process required to generate the model can be highly time consuming and somewhat challenging particularly when creating a model, or a digital twin of a physical test setup rather than relying solely on simulation data. In practice, the development of any reduced-order thermal model generally follows three stages, illustrated in Fig. \ref{fig:ROMstages}.  Stage A involves constructing a full-system CFD model or preparing a physical test setup, Stage B focuses on generating the ROM by identifying model parameters from simulation or experimental data, analogous to the training stage in machine-learning workflows and Stage C consists of using the ROM to evaluate temperature responses under new transient operating conditions. For 3D CFD analysis or laboratory testing performed without ROM development, Stage B is omitted.
In its original form, construction of an LPLSP model required performing a sequence of parametric studies in simulation, where each heat source was driven independently under constant input. Temperatures at all monitored locations were then extracted and fitted to the governing LPLSP equations to determine the model parameters. This workflow is labor intensive and scales poorly with the number of heat sources. Moreover, in many cases the prototype already exists \cite{neel_ecce, exp1, exp2} and the ability to monitor and control the temperature of experimental setup is critical, especially in fabrication and operation of the device \cite{kumar1, micro, kumar2, kumar3}. For such systems, applying isolated, constant-input excitation to individual components is impractical or infeasible. To overcome these limitations, the present work introduces an \emph{ensemble parameter estimation} approach that identifies all LPLSP parameters simultaneously from a single transient dataset. The method accommodates fully transient inputs or even pseudo-random input signals, thereby eliminating the need for controlled parametric excitation. This substantially reduces the effort and computation time required in Stage B and enables a high degree of automation in ROM generation. The remainder of this paper presents the formulation of the proposed approach, its implementation using rank-reduced and two-stage estimation strategies, and its application to pure conduction, natural convection, forced convection problems, and a case study of heat transfer in a power inverter module.
%%%%%%%%%%%%%%%%%%%%%%%%
\section{Overview of CFD simulation setup}
\begin{figure}[t]
	\centering
	% Row 1: a and b side-by-side
	\begin{subfigure}[b]{0.48\textwidth}
		\centering
		\includegraphics[width=0.9\linewidth]{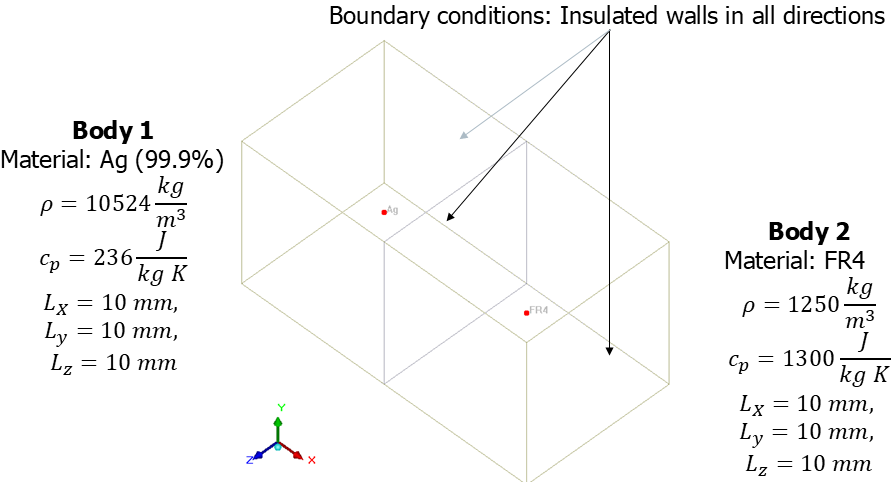}
		\caption{Setup for two-body conduction}
		\label{fig:case1}
	\end{subfigure}
	\hfill
	\begin{subfigure}[b]{0.48\textwidth}
		\centering
		\includegraphics[width=0.9\linewidth]{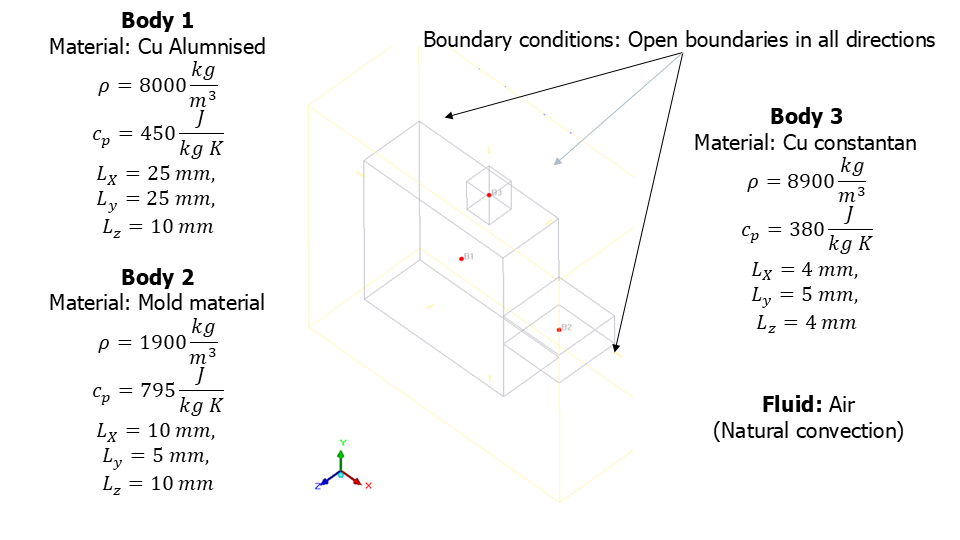}
		\caption{Setup for three-body convection}
		\label{fig:case2}
	\end{subfigure}
	
	% Row break
	\vspace{0.8em}
	
	% Row 2: c full width
	\begin{subfigure}[b]{0.98\textwidth}
		\centering
		\includegraphics[width=0.8\linewidth]{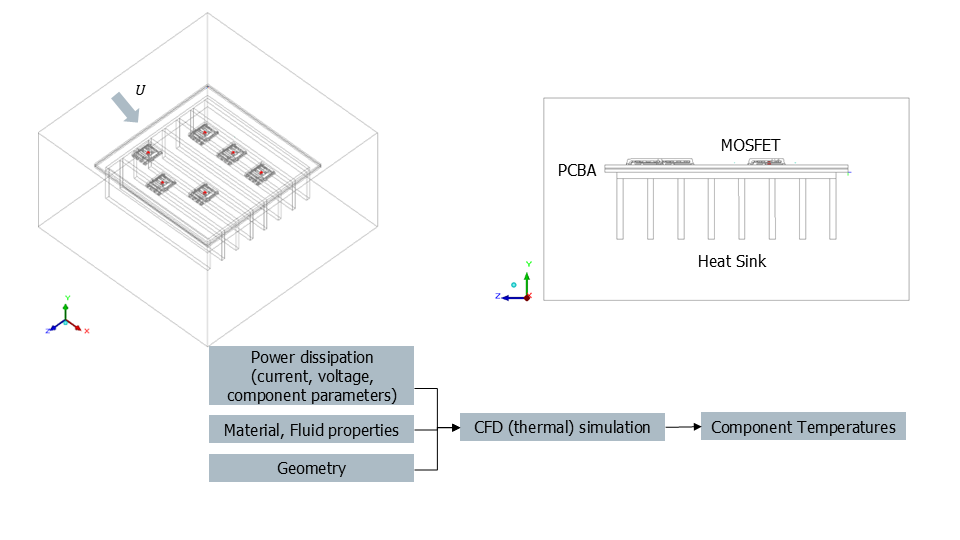}
		\caption{Setup for case study of inverter module with six MOSFETs on a PCBA attached to a heatsink.}
		\label{fig:case3}
	\end{subfigure}
	
	\caption{Overview of simulation setup}
	\label{fig:overview}
\end{figure}

The simulation setup described here is used to generate the full CFD data, training data for model development, and to calculate computational times. The case studies considered in this work include (i) a two-body conduction system, (ii) a three-body natural convection system, and (iii) an inverter module comprising six MOSFETs mounted on a printed circuit board assembly (PCBA) attached to a finned heatsink under natural convection and forced convection environment at constant flow velocity (Fig.~\ref{fig:overview}). These configurations represent a broad class of thermal management problems encountered in automotive and industrial electronics, ranging from simple conduction-dominated assemblies to complex multi-source systems with coupled conduction and convection. \\
All simulations were performed using Ansys\textsuperscript{\textregistered} Icepak\textsuperscript{\texttrademark}. The solver computes the incompressible Navier–Stokes and energy equations with buoyancy modeled via the Boussinesq approximation for natural convection, while applying ideal gas equation of state for forced convection. A zero-order turbulence model was employed, with first- and second-order upwind schemes for spatial discretization and a fully implicit scheme for temporal integration. Initial conditions were set to $20^{\circ}\mathrm{C}$ and $1~\mathrm{atm}$, with standard earth gravity. Open boundary conditions were applied for natural convection and inlet velocity of $U=10 m/s$ is applied for forced convection. Material properties for silicon, copper, FR4, and aluminum were assigned based on manufacturer data to ensure accurate thermal conduction paths. The material properties used in these simulations are identical to those reported in our previous work \cite{LPLSP2}. Mesh resolution was selected through a grid-independence study, with total cell counts ranging from $3\times10^5$ to $3\times10^6$ depending on geometry complexity. Grid stretching and clustering were applied near heat sources and critical flow regions, while non-conformal meshing allowed separate refinement of PCBA and heatsink assemblies. Cut-cell meshing was used to preserve geometric fidelity. These settings ensured convergence of flow and energy equations under the solver’s default criteria.

%%%%%%%%%%%%%%%%%%%%%%%%
\section{Ensemble parameter estimation for LPLSP model} % WORKING HERE
The Lumped Parameter Linear Superposition (LPLSP) method \cite{LPLSP1,LPLSP2} models transient temperatures by combining lumped thermal networks with the principle of linear superposition. Each monitored temperature is expressed as a sum of contributions from all heat sources, weighted by thermal resistances and time constants. For the cases of pure conduction, natural convection and forced convection with constant flow velocity, the relationship between the transient temperature of a body and the corresponding transient input power dissipation is defined by the LPLSP formulation as
\begin{equation}
	T_i(t) = T_0 + T_L (t) + \sum_{j=1}^{N} P_j(t)\, R_{ij}\,\bigl(1 - e^{-K_{ij} t}\bigr),
\end{equation}
where:
\begin{itemize}
	\item $T_0$ is the initial temperature,
	\item $T_L(t) = (P_T / C_T)t$ is the linear temperature component, with $P_T$ denoting the total input power from all sources and $C_T$ the total thermal capacitance of all objects (including working fluid),
	\item $N$ is the number of heat sources,
	\item $P_j(t)$ is the power dissipation of source $j$,
	\item $R_{ij}$ and $K_{ij}$ are the thermal resistance and time constant between monitor point $i$ and source $j$, respectively.
\end{itemize}
The linear temperature component depends on material properties and geometry and can be computed \textit{a priori}. For large systems with many sources and sinks, the total thermal capacitance becomes very large, making the contribution of the linear term negligible. For a two-body system, the conventional approach requires two separate simulations: one with $(P_1 = 1,\, P_2 = 0)$ and another with $(P_1 = 0,\, P_2 = 1)$. In the proposed method, instead of running individual parametric studies for each source, a single simulation is performed in which arbitrary transient or pseudo-random power dissipation inputs $P_1(t)$ and $P_2(t)$ are applied to the two bodies and the resulting temperatures $T_1(t)$ and $T_2(t)$ are recorded. Given the measured outputs $T_i(t)$, the inputs $P_j(t)$, and the model equation above, the only unknown parameters are the four thermal resistances $(R_{11}, R_{12}, R_{21}, R_{22})$ and the four time constants $(K_{11}, K_{12}, K_{21}, K_{22})$. These parameters are estimated using a nonlinear least-squares optimizer that minimizes the sum of squared residuals between the measured and modeled temperatures. In this work, the \texttt{SciPy.least-squares} function is used. Although the optimization method itself is identical to that used in the original LPLSP parameter-estimation process, the key difference is that all parameters are estimated simultaneously from a single simulation, as an ensemble estimation rather than through individual parametric runs. This results in a substantial reduction in computational effort and overall model development time. Once the parameters are identified, they are substituted back into the model equation, enabling rapid prediction of temperatures for new input power dissipation profiles. A pseudocode for this algorithm is presented in Algorithm \ref{pc1}. In this paper, in addition to presenting the results, the accuracy and computational effort of the proposed ensemble parameter estimation method are compared with those of the traditional approach, in which parameters are estimated individually for each source.
%%%
\begin{algorithm}[htpb]
	\caption{General $N$-Body Thermal Parameter Estimation - Naive implementation (without vectorization)}
	\begin{algorithmic}[1]
		% ------------------------------------------------------------
		\State \textbf{Input:} Measured temperatures (training dataset) with time $t$, power dissipations $P \in \mathbb{R}^{N \times T}$, measured temperatures $T_{\mathrm{meas}} \in \mathbb{R}^{M \times T}$, initial temperature $T_0 = 20^\circ$C
		\State Identify number of sources $N$ and temperature monitor points $M$ from the number of power/temperature channels
		% ------------------------------------------------------------
		\Statex
		\Function{ComputeTemperature}{$P, t, R_{i,:}, K_{i,:}, T_0$}
		%\Statex \hspace{1em}
		\State Initialize $T(t)=0$
		\For{$j = 1 \dots N$}                               
		\State Enforce $P_j(0)=0$.
		\State Identify indices where input changes $\mathcal{I} = \{ k : P_j(k) \ne P_j(k\!-\!1) \}$.
		%       \If{$\mathcal{I}$ is empty} \State \textbf{continue} \EndIf
		\State Compute temperature contribution of each source: $TC_j(k) = P_j(k)\,R_{ij}$.
		\For{each change index $k \in \mathcal{I}$}
		\State $\Delta TC = TC_j(k) - TC_j(k\!-\!1)$.
		\State $t_0 = t(k\!-\!1)$.
		\For{$m = k \dots T$}
		\State $T(m) = T(m) + \Delta TC \left(1 - e^{-K_{ij}(t(m)-t_0)}\right)$.
		\EndFor
		\EndFor
		\EndFor
		\State \Return $T_0 + T$
		\EndFunction
		% ------------------------------------------------------------
		\Statex
		\Function{Residuals}{$\theta$}
		\State $\theta$: parameter matrices: $R, K$
		%\State Initialize predicted temperature array $T_{\mathrm{pred}}$
		\For{$i = 1 \dots N$}
		\State Compute predicted temperature of node $i$:
		\[
		T_{\mathrm{pred},i}(t)
		= \Call{ThermalModel}{P, t, R_{i,:}, K_{i,:}, T_0}
		\]
		\EndFor
		\State Compute temperature errors: $E_i(t) = T_{\mathrm{pred},i}(t) - T_{\mathrm{meas},i}(t)$
		\State \Return Stack residuals $E_i(t)$ into a vector
		\EndFunction
		
		% ------------------------------------------------------------
		\Statex
		\State Choose initial guess $\theta_0$ for $R_{ij}$ and $K_{ij}$
		\State Solve nonlinear least-squares problem (\texttt{SciPy.least-squares}):
		\[
		\theta^\star = \arg\min_{\theta} \|\Call{Residuals}{\theta}\|_2^2
		\]
		\Statex
		\State \textbf{Output:} Estimated matrices $R^{\star}, K^{\star}$.
	\end{algorithmic}
	\label{pc1}
\end{algorithm}
%%%%%%%%%%%%%%%%%%
\subsection{Practical challenges in large thermal systems}
Although the proposed model is computationally efficient for systems with a small number of heat sources, its scalability becomes challenging as the system size increases. For a system with $N$ bodies, the total number of parameters to be identified grows as $2N^2$, because of two sets of parameters $R_{ij}, K_{ij}$. In simple configurations, the time constants can be approximated as $K_{ij} \approx 1/R_{ij} C_T$, where $C_T = \sum_i ^N \rho_i C_{p,i} V_i$ is the total thermal capacitance of the system. However, for more complex systems this approximation becomes unreliable, as the model is highly sensitive to the accuracy of the time constants. Consequently, both $R_{ij}, K_{ij}$ must be estimated independently rather than inferred from one another. In real-time applications such as the power inverter module with six heat sources, this requirement leads to a total of 72 parameters that must be estimated, significantly increasing the complexity of the identification problem. From extensive testing, it was observed that estimating all 72 parameters through a direct brute-force approach resulted in prohibitively long computation times ($t \sim 10^3 s$) and, in many cases, the solver still failed to reach full convergence. To address these limitations, a few improvements were incorporated: 
\begin{itemize}
	\item Vectorization of the core computations and acceleration using \texttt{numba} just-in-time compiler,
	\item Enforcing symmetry in the $R$ and $K$ parameter matrices.
\end{itemize}
These modifications substantially improved computational efficiency and convergence of solution.
%%%%%%%%
\subsection{Computational considerations and acceleration methods}
Some of the strategies used to accelerate the parameter estimation procedure are presented in this section. First, vectorization was applied to computation of model temperatures as functions of the transient input power dissipation. Although a naïve implementation based on nested loops is straightforward, its computational cost is prohibitively high. Reformulating the problem as a sequence of matrix-vector operations enables efficient vectorization, resulting in a substantial reduction in code execution time. Additional gains were achieved through the use of \texttt{numba} for just-in-time compilation of the core numerical routines. The most significant speed-up, however, arises from analyzing the structure of the $R$ and $K$ parameter matrices. From the two-body and three-body cases, it was observed that the resistance matrix $R_{ij} \in \mathbb{R}^{N \times N}$ is approximately symmetric about its main diagonal ($R_{ij} \approx R_{ji}$). Physically, this reflects the fact that the thermal resistance between a heat source and a measurement node is nearly identical regardless of the direction of heat transfer. A similar symmetry was observed in the time-constant matrix $K_{ij}$. Exploiting this property allows us to estimate only the upper triangular elements of each matrix and impose symmetry. This reduces the number of unique parameters per matrix from $N^2$ to $N(N+1)/2$. For a system with six heat sources and six temperature nodes, the total number of parameters therefore decreases from $2N^2 = 72$ to $2\,N(N+1)/2 = 42$. By combining vectorization, \texttt{numba} acceleration, and symmetry exploitation, the overall computation time is reduced from $t \sim 10^3 s$ to approximately $t \sim 10^1 s$. A pseudocode for this algorithm is presented in Algorithm \ref{pc2}.
\begin{algorithm}[htpb]
	\caption{$N$-Body Parameter Estimation by application of symmetry in $R,K$ and vectorization of core computations}
	\begin{algorithmic}[1]
		% ------------------------------------------------------------
		\State \textbf{Input:} Measured temperatures (training dataset) with time $t$, power dissipations $P \in \mathbb{R}^{N \times T}$, measured temperatures $T_{\mathrm{meas}} \in \mathbb{R}^{M \times T}$, initial temperature $T_0 = 20^\circ$C
		\State Identify number of sources $N$ and temperature monitor points $M$ from the number of power/temperature channels. In this case $M=N$ where the monitor points are co-located on the source.
		% ------------------------------------------------------------
		\Statex
		\Function{ComputeTemperatureVect}{$P, t, R_{i,:}, K_{i,:}, T_0$} (\textit{Vectorized core computations})
		%\Statex \hspace{1em}
		\State Initialize $T(t) = 0$
		\For{$j = 1 \dots N$}
		\State Enforce $P_j(0)=0$
		\State Identify indices where input changes: $\mathcal{I} = \{ k : P_j(k) \ne P_j(k-1) \}$
		\State Compute temperature contribution of each source: $TC_j(k) = P_j(k)\,R_{ij}$
		\For{each change index $k \in \mathcal{I}$}
		\State $\Delta TC = TC_j(k) - TC_j(k-1)$
		\State $t_0 = t(k-1)$
		\State Vectorized time differences: $\Delta t(m) = t(m) - t_0,\qquad m = k \dots T$
		\State Vectorized exponential kernel: $K_{\exp}(m) = 1 - e^{-K_{ij}\,\Delta t(m)}$
		\State Vectorized accumulation: $T(m) \gets T(m) + \Delta TC \cdot K_{\exp}(m),\quad m = k \dots T$
		\EndFor
		\EndFor
		\State \Return $T_0 + T$
		\EndFunction
		%-----------------------------------------------------------
		\Statex
		\Function{SymmetricParam}{$\theta, N$}
		\State Initialize $R, K \in \mathbb{R}^{N \times N}$.
		\State $n_{\mathrm{half}} = N(N+1)/2$.
		\State Extract vectors $R_{\mathrm{flat}}, K_{\mathrm{flat}}$ (upper triangles).
		\State Use only upper triangle and main diagonal: $R_{ij}=R_{\mathrm{flat}}(k), K_{ij}=K_{\mathrm{flat}}(k); i \le j $
		\State Apply symmetry about main diagonal: $R_{ji} = R_{ij}$, $K_{ji} = K_{ij}$.
		\State \Return $R,K$.
		\EndFunction
		% -----------------------------------------------------------
		\Statex
		\Function{Residuals}{$\theta$}
		\State $(R,K) = \Call{SymmetricParam}{\theta,N}$.
		\For{$i=1 \dots N$}
		\[
		T_{pred_{i}}(t) = \Call{ComputeTemperature}{P,t,R_{i,:},K_{i,:}, T_0}
		\]
		\EndFor
		\State \Return $\mathrm{vec}(T_{pred} - T_{meas})$.
		\EndFunction
		% ------------------------------------------------------------
		\Statex
		\State Choose initial guesses:  $\theta_0 = \mathbf{1}_{\,2 n_{\mathrm{half}} \times 1}$
		\State Solve non-linear least squares problem (\texttt{Scipy.least-squares}):
		\[
		\theta^\star = \arg\min_{\theta} \|\Call{Residuals}{\theta}\|_2^2
		\]
		\State Symmetric matrices: $(R^\ast, K^\ast) = \Call{SymmetricParam}{\theta^\ast, N}.$
		%------------------------------------------------------------
		\Statex
		\State \textbf{Output:} Estimated symmetric matrices $R^{\star},K^{\star} \in \mathbb{R}^{N \times N}$
	\end{algorithmic} \label{pc2}
\end{algorithm}
%%%%%%%%
\subsection{Parameter estimation for rectangular coupling matrices}
The acceleration strategies described in the previous section are applicable primarily to square $R$ and $K$ matrices, which arise when temperature monitoring points are located directly on the heat sources or in very close proximity. In such configurations, there is strong thermal coupling between the monitor points of different sources, and the thermal resistance between two sources is approximately symmetric. This justifies the assumption $R_{ij} = R_{ji}$ and $K_{ij} = K_{ji}$. However, in cases where the monitor points are not located on the sources or placed farther away on a heatsink or printed circuit board assembly (PCBA), the number of monitor points $M$ may exceed the number of sources $N$, resulting in rectangular $R$ and $K$ matrices of size $N \times M$. For such rectangular matrices, parameter reduction via symmetry about the main diagonal is no longer valid. While it may be possible to artificially pad the source with zero vectors to match $M=N$ and create a square matrix, this approach is physically inconsistent. In pure conduction systems, where heat is not removed from the system through other mechanisms, the thermal resistance between a heat source and a sink can be assumed to be reciprocal. However, when heat is removed through mechanisms such as convection or radiation, reciprocity no longer holds ($R_{ij} \neq R_{ji}, K_{ij} \neq K_{ji}$). In such cases, the sinks introduce directional behavior with non-reciprocal boundary conditions, since convection and radiation are inherently unidirectional. \\
For the inverter module with 6 sources and 8 monitor points (including the sources, PCBA, and heatsink), the $R$ and $K$ matrices become rectangular, with 48 elements each, resulting in a total of 96 parameters. In this case, even the combination of vectorization and \texttt{numba} acceleration is insufficient to achieve substantial reductions in computation time, as the number of parameters increases from 42 (for the square case) to 96. Therefore, when dealing with rectangular matrices, two additional strategies are employed to reduce computation time.
\subsubsection{Rank reduction of $R$ and $K$ matrices}
This approach is conceptually similar to Proper Orthogonal Decomposition (POD) \cite{POD3} or Principal Component Analysis (PCA), where a high-dimensional system is efficiently represented using a small number of dominant modes. In this context, the matrices $R$ and $K$ act as spatial coupling operators between the heat sources and the temperature monitoring points. A low-rank approximation is physically appropriate because, in complex PCB or heatsink geometries, many temperature nodes exhibit similar spatial behavior, and multiple heat sources influence monitor points in correlated patterns. Therefore, only a few spatial basis functions are typically sufficient to represent all rows of $R$, and the same holds for $K$, since spatial variations in time constants tend to evolve smoothly across the domain. In essence, the temperature fields vary smoothly with position, the thermal influence of each source decays gradually with distance, and many rows of $R$ and $K$ are strongly correlated. Thus, most of the system’s energy can be captured by a small number of dominant modes. To exploit this structure, the rectangular matrices are approximated using low-rank factorizations, $R = AB^{T}$ and $K = CD^{T}$, where $A \in \mathbb{R}^{M \times r}$ and $B \in \mathbb{R}^{N \times r}$, where $r$ denotes the chosen rank. For the inverter module with six heat sources ($N$) and eight monitoring points ($M$), the full matrices contain 48 parameters each. A rank-3 representation reduces this to 42 parameters $(A: 24, B: 18)$, and a rank-2 representation reduces it further to only 28 parameters. This reduction dramatically accelerates optimization and results in smoother estimates while accommodating the rectangular structure of the matrices. From the case studies it is observed that even a rank-1 approximation works well, rank-2 is almost always sufficient, and rank-3 is nearly indistinguishable from the full model. In the reduced formulation, each element of the coupling matrices is expressed as, $R_{ij} = \sum_{k=1}^{r} A_{ik} B_{jk}, K_{ij} = \sum_{k=1}^{r} C_{ik} D_{jk},$ and only the elements of $A$, $B$, $C$, and $D$ are optimized. A pseudocode for the rank-reduced method is presented in Algorithm \ref{pc3}. The rank estimation can also be done algorithmically by computing the singular value decomposition (SVD) of $R, K$ matrices. The singular values $\sigma_i$, indicate the relative contribution of each mode, and enables selection of the smallest number of modes that captures a desired fraction $\tau$ of the total information. 
\begin{equation}
	\begin{split}
		r^\star &= \min \left\{ k \;\middle|\;
		\sum_{i=1}^{k} \frac{\sigma_i}{\sum_{j} \sigma_j} \ge \tau \right\}, \\
		R &= \sum_{i=1}^{r^\star} \sigma_i \, \mathbf{u}_i \mathbf{v}_i^{\top}, 
		\qquad
		\sigma_1 \ge \sigma_2 \ge \cdots \ge \sigma_{r^\star} .
	\end{split}
\end{equation}

\begin{algorithm}[htpb]
	\caption{Rank-Reduced Parameter Estimation for system with $N$-Inputs, $M$-Temperature monitor points}
	\begin{algorithmic}[1]
		\State \textbf{Input:} Measured temperatures (training dataset) with time $t$, power dissipations $P \in \mathbb{R}^{N \times T}$, measured temperatures $T_{\mathrm{meas}} \in \mathbb{R}^{M \times T}$, initial temperature $T_0 = 20^\circ$C
		%\vspace{1mm}
		\Statex \textbf{Compute Temperature:}  Function \textsc{ComputeTemperature}$(P,t,R_{i,:},K_{i,:},T_0)$ 
		defined previously in Algorithm \ref{pc2}.
		%\vspace{1mm}
		\Statex \textbf{Rank-Decomposition of Parameters:}
		\State Choose initial rank $r$
		\Function{ReduceRank}{$\theta$}
		\State Reshape $\theta$ into matrices $A \in \mathbb{R}^{M \times r}$, $B \in \mathbb{R}^{N \times r}$
		\State Reshape remaining parameters into $C \in \mathbb{R}^{M \times r}$, $D \in \mathbb{R}^{N \times r}$
		\State Form $R = A B^{T}$ and $K = C D^{T}$
		\State \Return $R,K$
		\EndFunction
		%\vspace{1mm}
		\Statex \textbf{Residual Computation:}
		\Function{Residuals}{$\theta$}
		\State $(R,K)=\Call{ReduceRank}{\theta}$
		\For{$i=1\dots M$}
		\State $T^{\text{pred}}_i(t)=\Call{ComputeTemperature}{P,t,R_{i,:},K_{i,:},T_0}$
		\EndFor
		\State \Return $\mathrm{vec}(T_{\mathrm{pred}} - T_{\mathrm{meas}})$
		\EndFunction
		%\vspace{1mm}
		\State Choose initial guess $\theta_0$
		\State Solve nonlinear least squares:
		\[
		\theta^\star = \arg\min_{\theta} \|\Call{Residuals}{\theta}\|_2^2
		\]
		\State $(R^\star,K^\star) = \Call{ReduceRank}{\theta^\star}$
		%\vspace{1mm}
		\State \textbf{Output:} Estimated low-rank matrices $R^\star$ and $K^\star$
	\end{algorithmic} \label{pc3}
\end{algorithm}

\subsubsection{Two-stage parameter estimation}
The second approach employs a two-stage parameter estimation process. For the inverter module, Stage 1 focuses only on the square subsystem comprising the six heat sources and their corresponding (co-located) six temperature monitoring points. In this stage, we estimate the $6 \times 6$ matrices ($R$ and $K$), where conduction is the dominant mechanism. Stage 1 utilizes only the datasets $(T_1\text{--}T_6)$ and $(P_1\text{--}P_6)$. For this subset of the system, we impose symmetry about the main diagonal and estimate only the parameters in the upper triangular portions of the $R$ and $K$ matrices. These parameters are then fixed and subsequently used to estimate the thermal couplings to the remaining sinks (the PCBA and the heatsink). In Stage 2, with $R(1{:}6,1{:}6)$ and $K(1{:}6,1{:}6)$ held constant, we estimate only the couplings for the PCBA and heatsink. The temperature responses for the additional monitor points, $T_7(t)$ and $T_8(t)$, are modeled as,
\begin{equation}
	\begin{split}
		T_7(t) = T_0 + \sum_{j=1}^{6} P_j(t)\, R_{7j} \left(1 - e^{-K_{7j} t}\right), \\
		T_8(t) = T_0 + \sum_{j=1}^{6} P_j(t)\, R_{8j} \left(1 - e^{-K_{8j} t}\right).
	\end{split}
\end{equation}
This approach is physically consistent, highly accurate and takes significantly low computation time. A pseudocode for this is presented in Algorithm \ref{pc4}. In addition to this, we can also combine both these approaches of two stage estimation and rank reduction to speed up the process for very large systems. 
%%%%%%
\begin{algorithm}[htpb]
	\caption{Two-Stage Parameter Estimation for Inverter module with rectangular $R,K$ matrices ($N=6$ sources, $M=8$ monitor points)}
	\begin{algorithmic}[1]
		% ------------------------------------------------------------
		\State \textbf{Input:} 
		Measured data $t$, power inputs $P \in \mathbb{R}^{6 \times T}$, temperature measurements $T_{\mathrm{meas}} \in \mathbb{R}^{8 \times T}$, initial temperature $T_0 = 20^\circ$C.
		\State Split temperature data:
		\[
		T_{\mathrm{src}} = T_{\mathrm{meas}}(1{:}6,:), \qquad
		T_{\mathrm{sink}} = T_{\mathrm{meas}}(7{:}8,:).
		\]
		%\vspace{1mm}
		\Statex \textbf{Compute Temperature:}  Function \textsc{ComputeTemperature}$(P,t,R_{i,:},K_{i,:},T_0)$ 
		defined previously in Algorithm \ref{pc2}.
		\vspace{1mm}
		\Statex \textbf{Symmetric Parameter}  Function \textsc{SymmetricParam}$(\theta,N)$ defined previously in Algorithm \ref{pc2}.
		%\vspace{1mm}
		\Statex \textbf{Residuals for Stage 1: Estimate symmetric source block}
		\Function{ResidualsStage1}{$\theta,t,P_{\mathrm{src}},T_{\mathrm{src}},T_0$}
		\State $(R,K)=\Call{SymmetricParam}{\theta,6}$
		\For{$i=1\dots 6$}
		\State 
		$T^{\mathrm{pred}}_i(t)=
		\Call{ComputeTemperature}{P_{\mathrm{src}},t,R_i,K_i,T_0}$
		\EndFor
		\State \Return $\mathrm{vec}(T_{\mathrm{pred}} - T_{\mathrm{src}})$
		\EndFunction
		%\vspace{1mm}
		\State Choose initial guess $\theta_0$ for the symmetric block.
		\State Solve:
		\[
		\theta_1^\star 
		= 
		\arg\min_{\theta}\|\Call{ResidualsStage1}{\theta}\|_2^2.
		\]
		\State $(R_{\mathrm{src}},K_{\mathrm{src}}) 
		= \Call{SymmetricParam}{\theta_1^\star,6}$
		%\vspace{1mm}
		\Statex \textbf{Residuals for Stage 2: Couplings to Sink Nodes}
		\Function{ResidualsStage2}{$\theta,t,P_{\mathrm{src}},T_{\mathrm{sink}},T_0,R_{\mathrm{src}},K_{\mathrm{src}}$}
		\State Reshape $\theta$ into
		\[
		R_{\mathrm{sink}} \in \mathbb{R}^{2 \times 6}, 
		\qquad 
		K_{\mathrm{sink}} \in \mathbb{R}^{2 \times 6}
		\]
		\For{$i = 1,2$} 
		\State
		$T^{\mathrm{pred}}_i(t) 
		= \Call{ComputeTemperature}{P_{\mathrm{src}},t,
			R_{\mathrm{sink}}(i,:),K_{\mathrm{sink}}(i,:),T_0}$
		\EndFor
		\State \Return $\mathrm{vec}(T_{\mathrm{pred}} - T_{\mathrm{sink}})$
		\EndFunction
		%\vspace{1mm}
		\State Choose initial guess $\theta_0^{(2)}$ for sinks.
		\State Solve:
		\[
		\theta_2^\star 
		= 
		\arg\min_{\theta}\|\Call{ResidualsStage2}{\theta}\|_2^2.
		\]
		\State Extract $R_{\mathrm{sink}},K_{\mathrm{sink}}$ from $\theta_2^\star$
		%\vspace{1mm}
		\State \textbf{Output:}
		\[
		R^\star =
		\begin{bmatrix}
			R_{\mathrm{src}} \\[2mm]
			R_{\mathrm{sink}}
		\end{bmatrix},
		\qquad
		K^\star =
		\begin{bmatrix}
			K_{\mathrm{src}} \\[2mm]
			K_{\mathrm{sink}}
		\end{bmatrix}.
		\]
	\end{algorithmic}\label{pc4}
\end{algorithm}
%%%%%%%%%%%%%%
\section{Results}
\begin{figure}[htbp]
	\centering
	
	% --- Row 1 ---
	\begin{subfigure}[b]{0.48\textwidth}
		\centering
		\includegraphics[width=\textwidth]{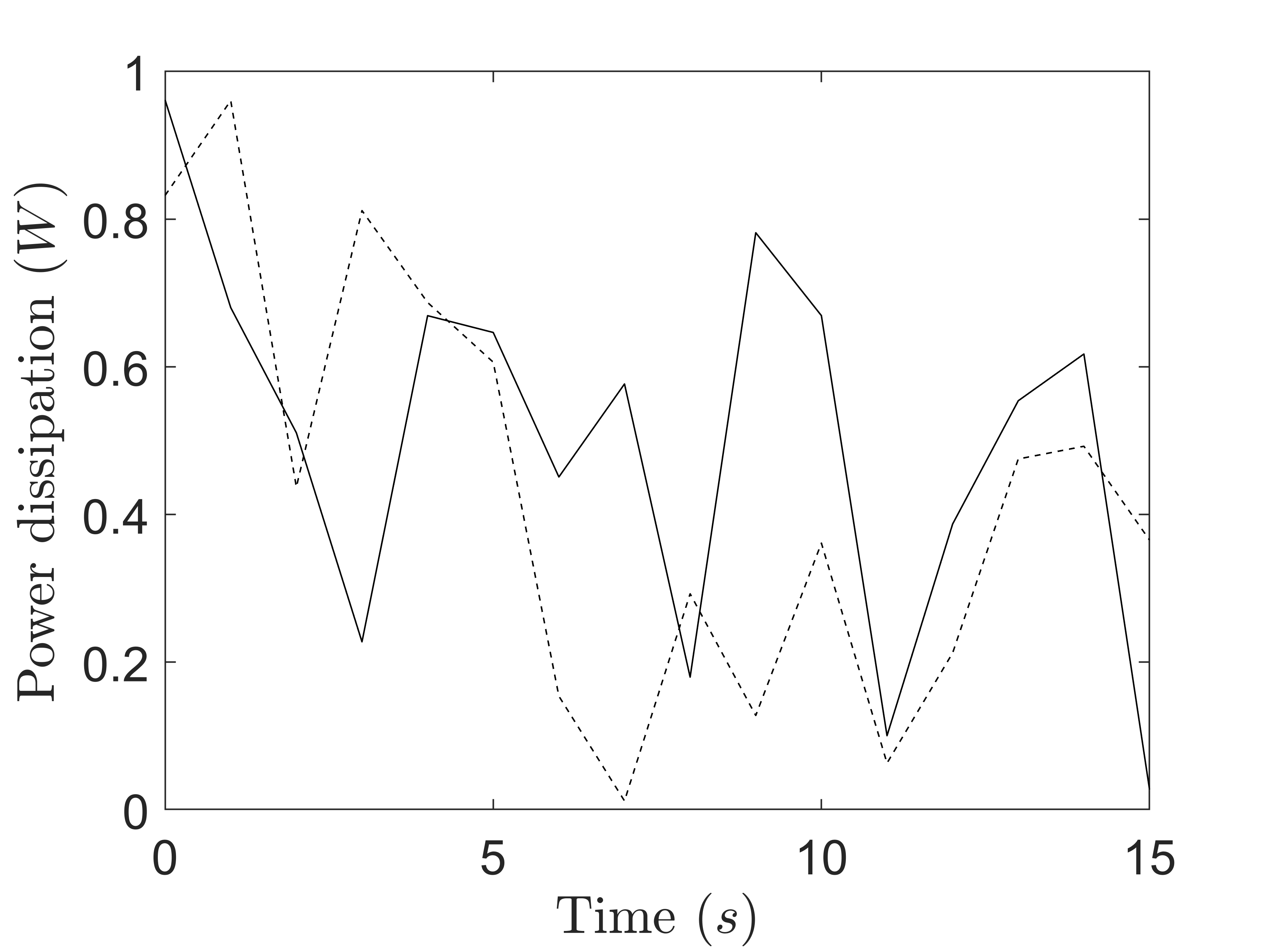}
		\caption{Model development: Power dissipation $P_1 (-), P_2 (--)$}
		\label{fig:Case1_Training_P}
	\end{subfigure}
	\hfill
	\begin{subfigure}[b]{0.48\textwidth}
		\centering
		\includegraphics[width=\textwidth]{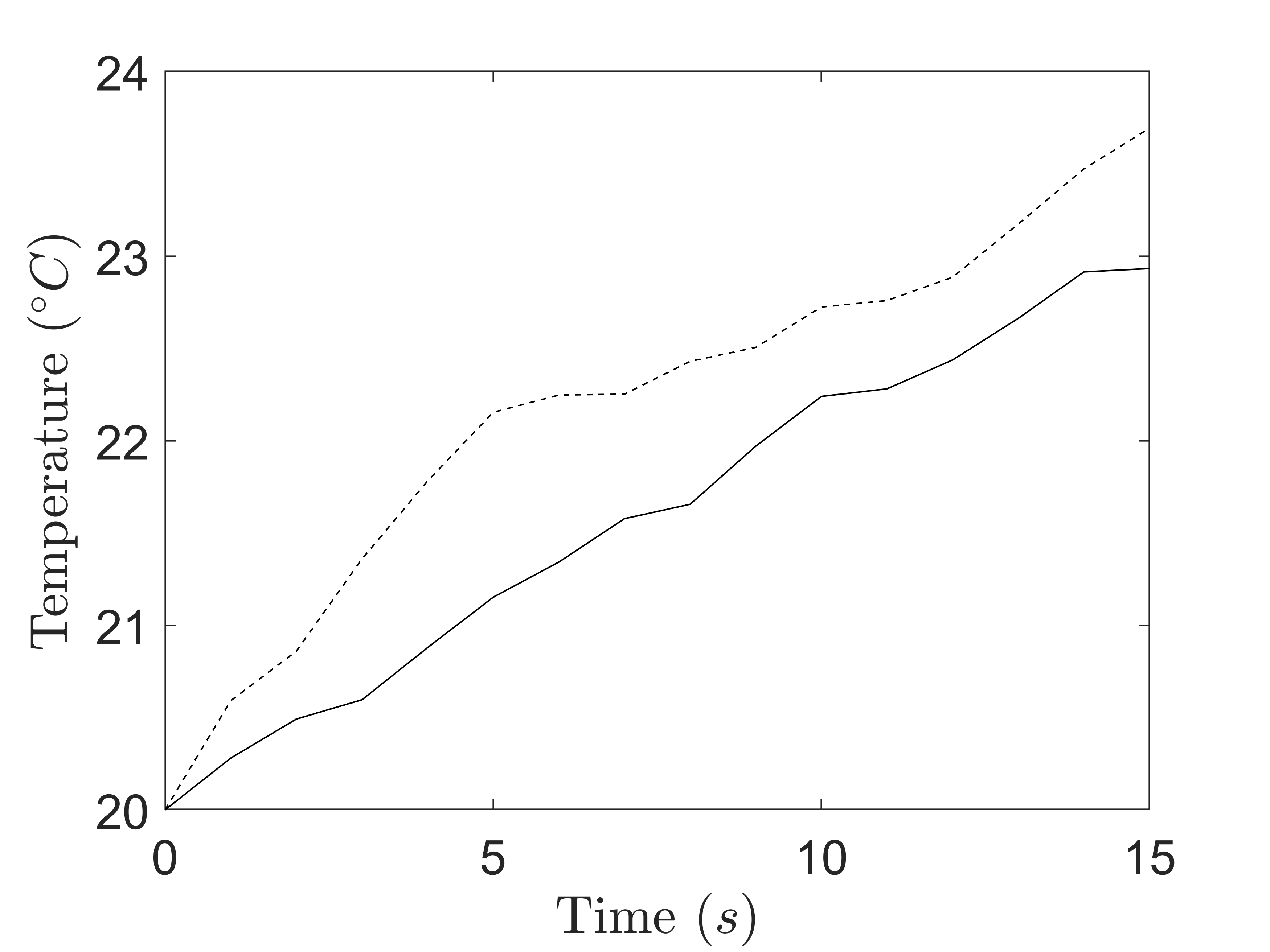}
		\caption{Model development: Temperature output $T_1 (-), T_2 (--)$}
		\label{fig:Case1_Training_T}
	\end{subfigure}
	\par\medskip % safer than \\ inside figure
	% --- Row 2 ---
	\begin{subfigure}[b]{0.48\textwidth}
		\centering
		\includegraphics[width=\textwidth]{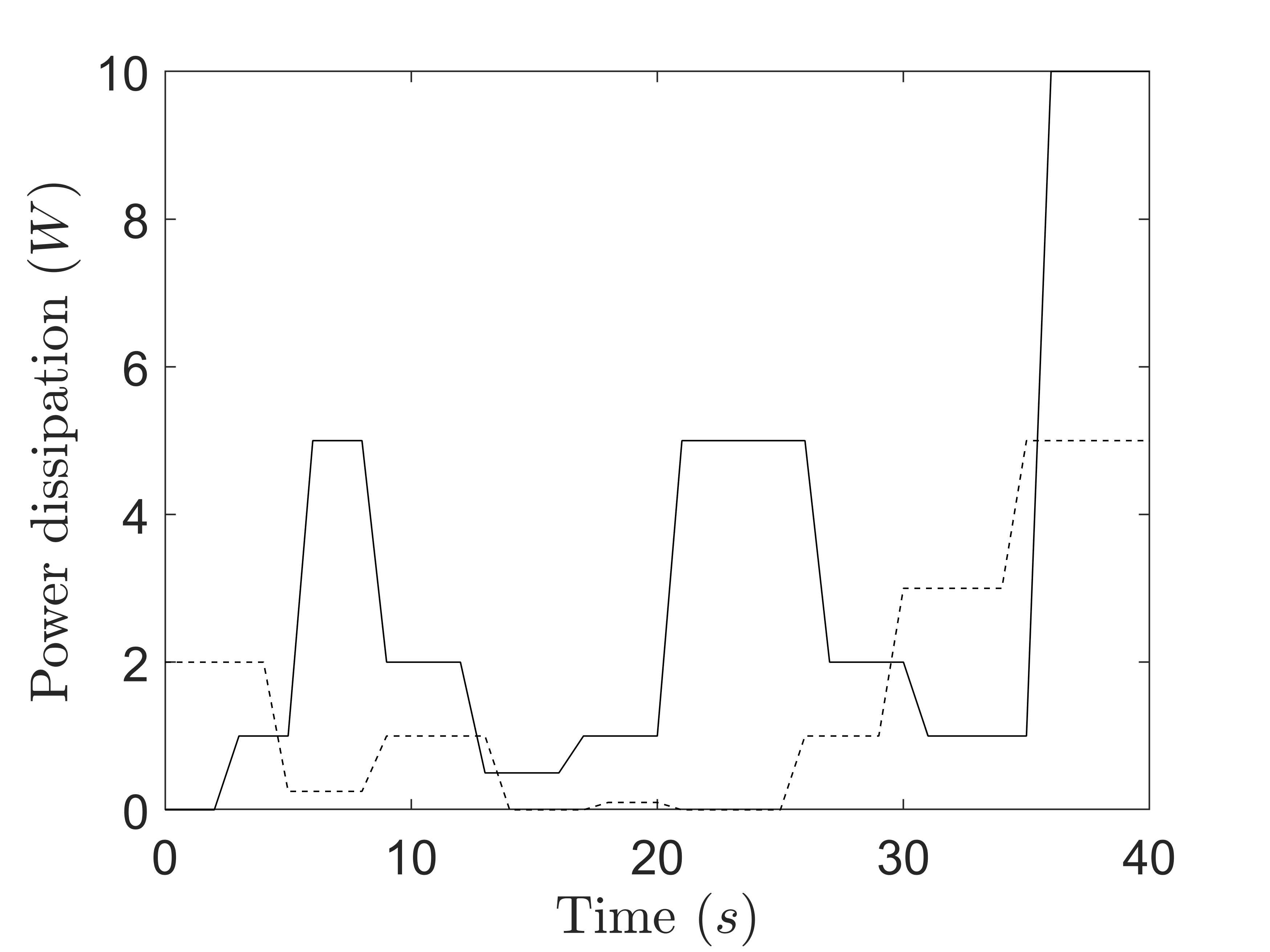}
		\caption{Validation: Applied power dissipation $P_1 (-), P_2 (--)$}
		\label{fig:Case1_Test_P}
	\end{subfigure}
	\hfill
	\begin{subfigure}[b]{0.48\textwidth}
		\centering
		\includegraphics[width=\textwidth]{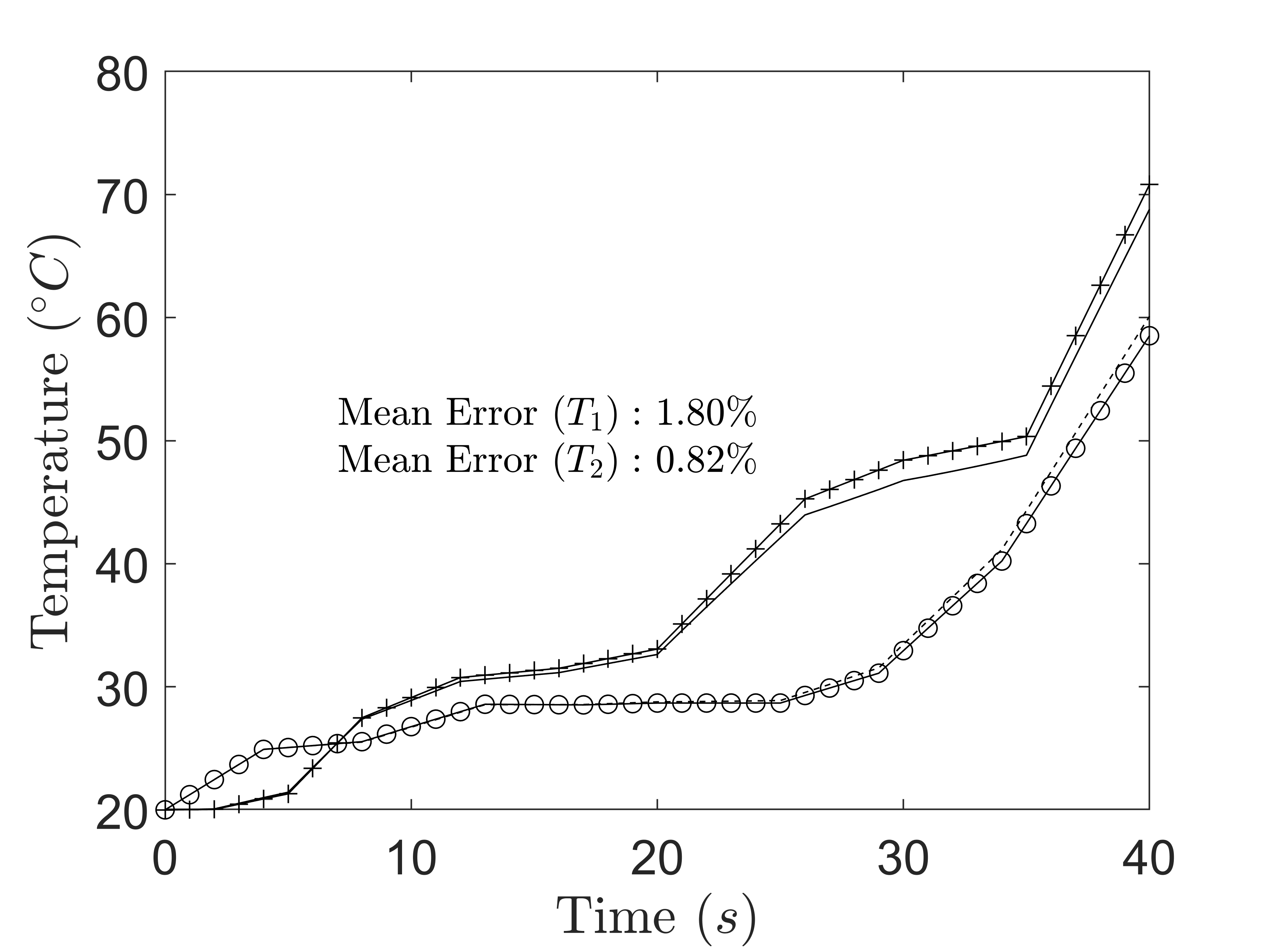}
		\caption{Validation: Simulation temperature $T_1 (-), T_2 (--)$ vs model $T_1 (-+), T_2 (-o)$}
		\label{fig:Case1_Test_T}
	\end{subfigure}
	\caption{Comparison of temperatures from simulation and LPLSP model for the two-body conduction case. The generation of training / model-development data via CFD requires $162 s$. A traditional parametric study would incur this cost multiplied by the number of heat sources (2 in this case). The mean percentage error between the model predictions and the simulation results is also reported.}
	\label{fig:Case1}
\end{figure}

\begin{figure}[htbp]
	\centering
	% --- Row 1 ---
	\begin{subfigure}[b]{0.48\textwidth}
		\centering
		\includegraphics[width=\textwidth]{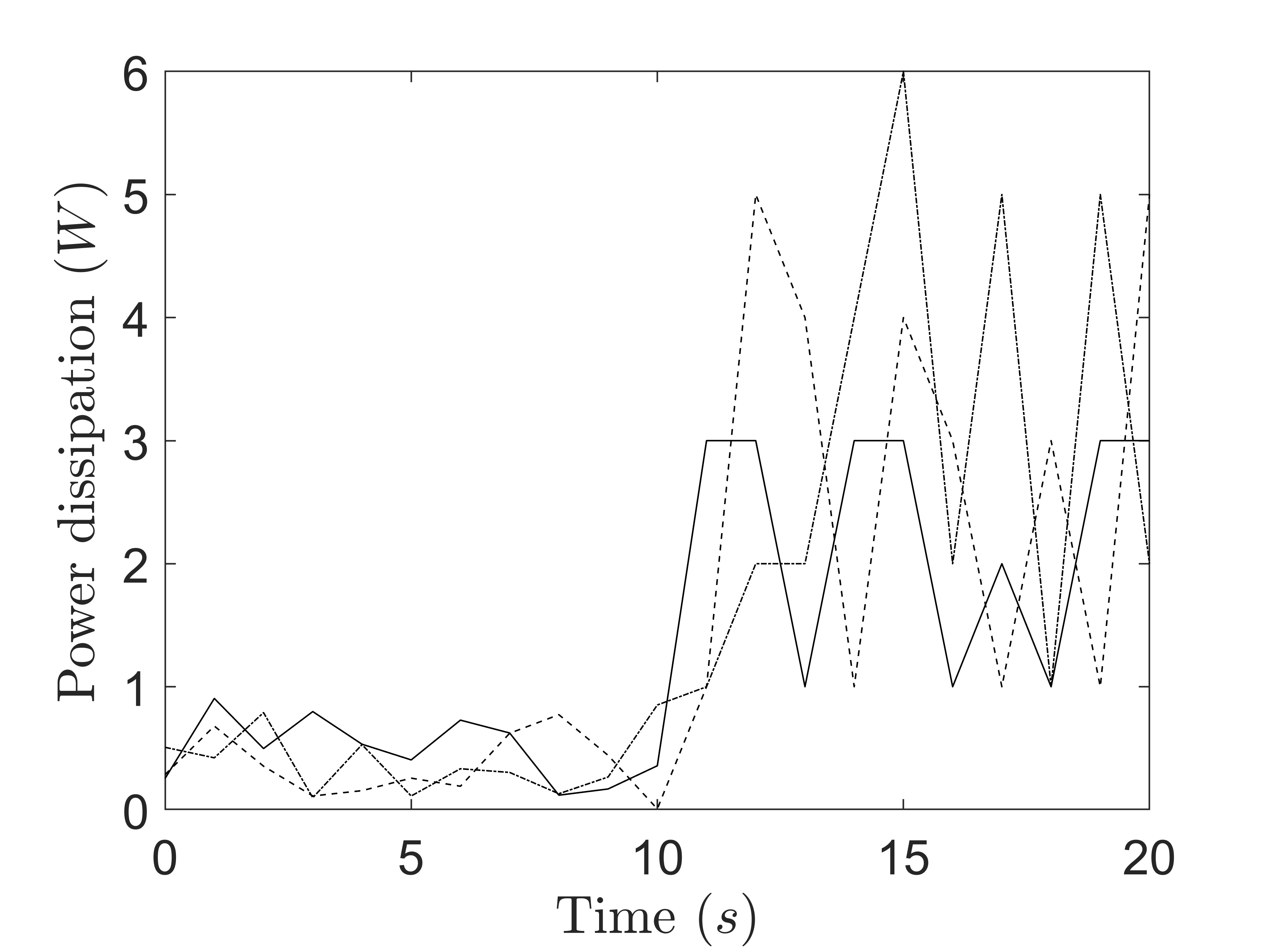}
		\caption{Model development: Input $P_1 (-), P_2 (--), P_3 (-.)$}
		\label{fig:Case2_Training_P}
	\end{subfigure}
	\hfill
	\begin{subfigure}[b]{0.48\textwidth}
		\centering
		\includegraphics[width=\textwidth]{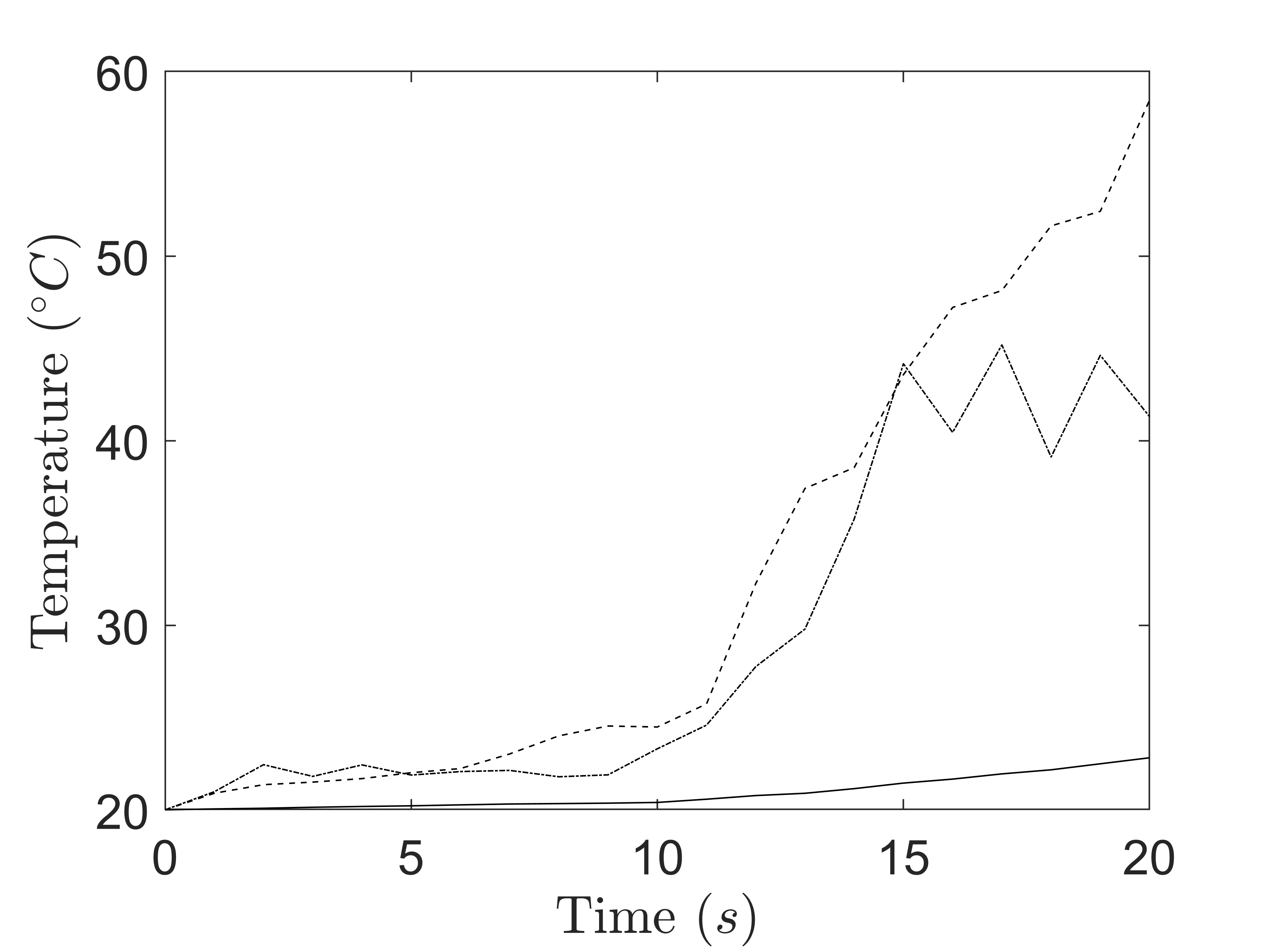}
		\caption{Model development: Output $T_1 (-), T_2 (--), T_3 (-.)$}
		\label{fig:Case2_Training_T}
	\end{subfigure}
	\par\medskip % safer than \\ inside figure
	% --- Row 2 ---
	\begin{subfigure}[b]{0.48\textwidth}
		\centering
		\includegraphics[width=\textwidth]{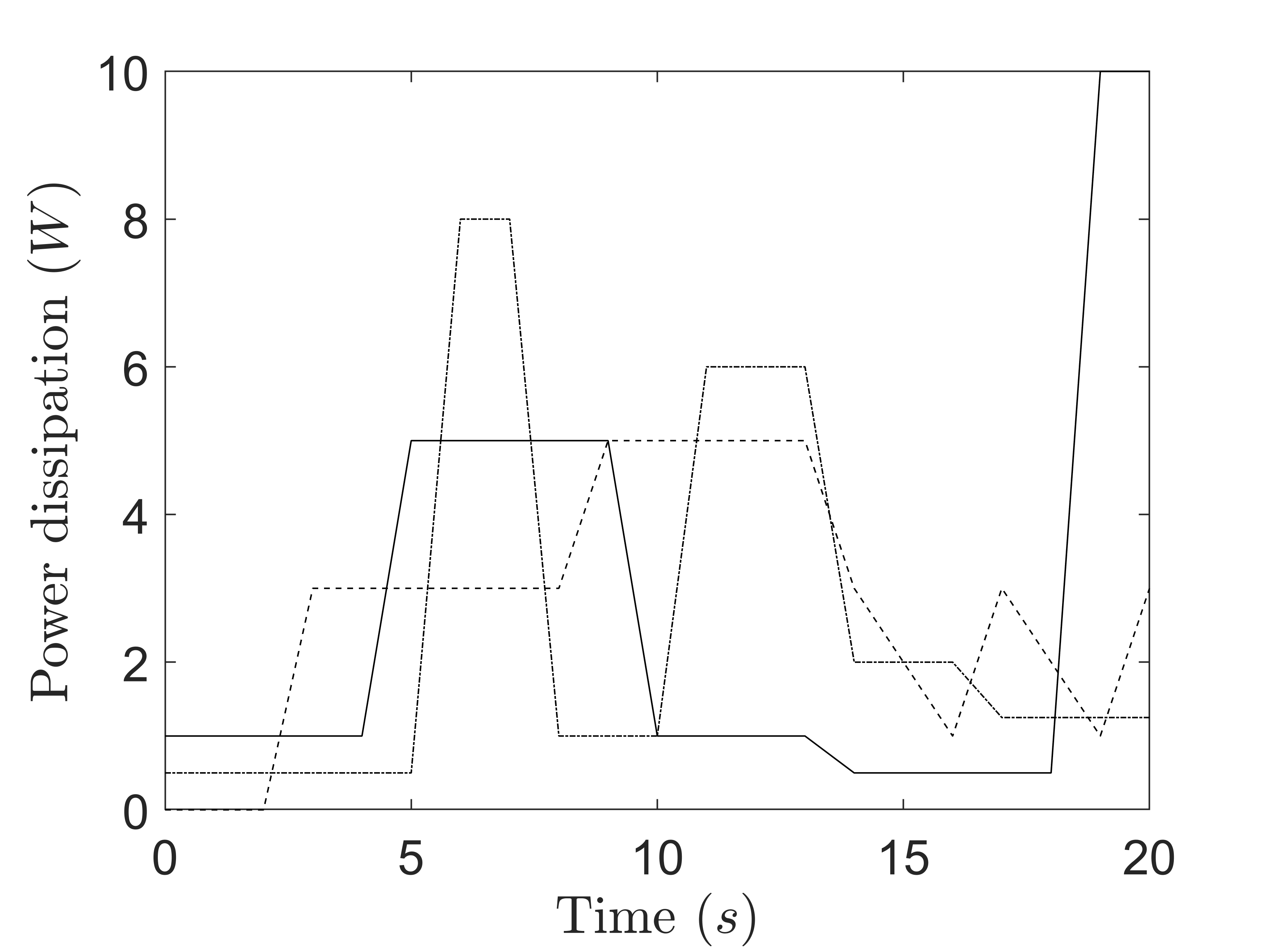}
		\caption{Validation: Input $P_1 (-), P_2 (--), P_3 (-.)$}
		\label{fig:Case2_Test_P}
	\end{subfigure}
	\hfill
	\begin{subfigure}[b]{0.48\textwidth}
		\centering
		\includegraphics[width=\textwidth]{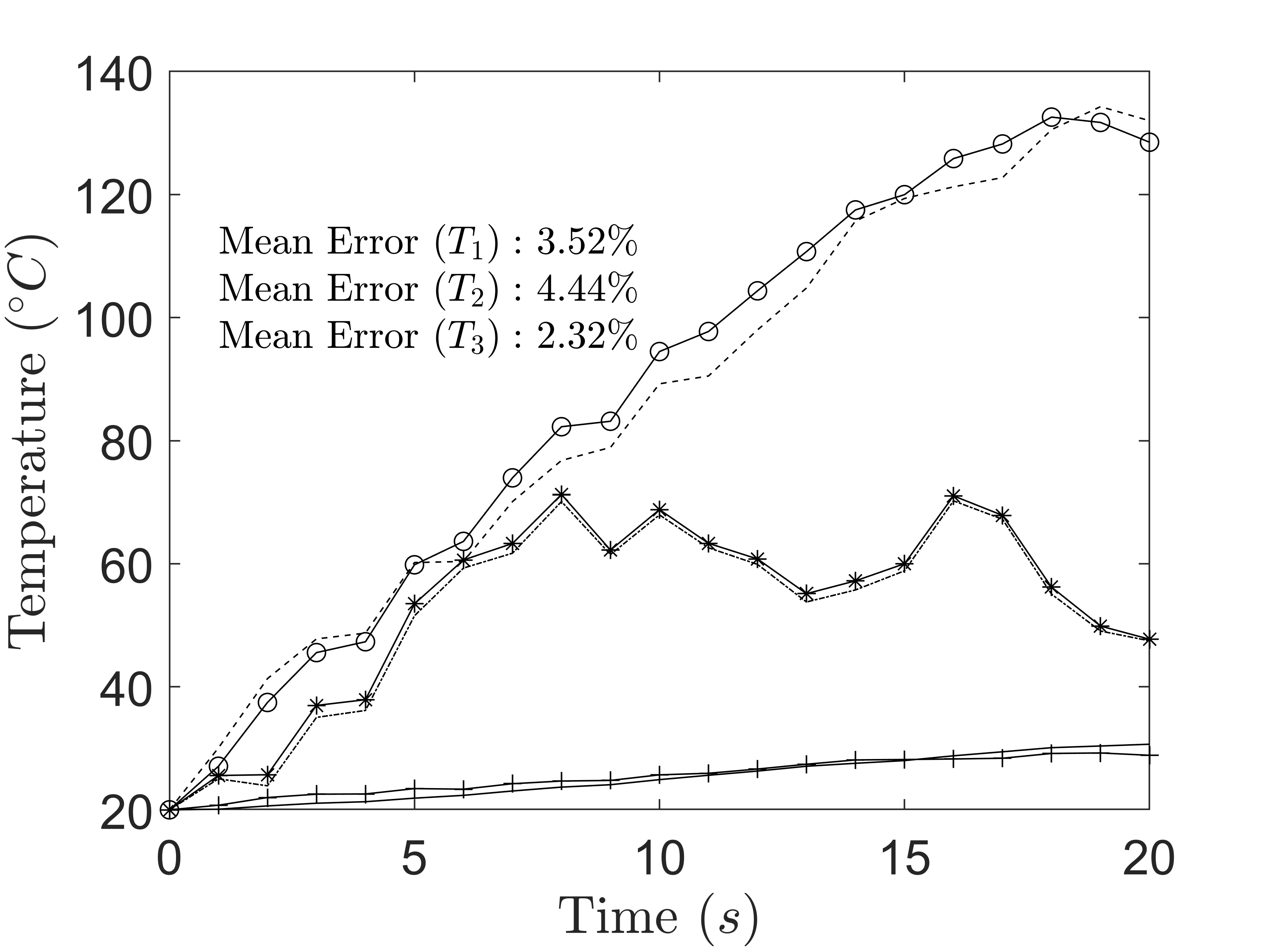}
		\caption{Validation: Simulation temperature $T_1 (-), T_2 (--) , \newline T_3 (-.)$ vs model $T_1 (-+), T_2 (-o), T_3 (-*)$}
		\label{fig:Case2_Test_T}
	\end{subfigure}
	\caption{Comparison of temperatures from simulation and LPLSP model for the two-body conduction case. The generation of training / model-development data via CFD requires $425 s$. A traditional parametric study would incur this cost multiplied by the number of heat sources (3 in this case). The mean percentage error between the model predictions and the simulation results is also reported.}
	\label{fig:Case2}
\end{figure}

\begin{figure}[htbp]
	\centering
	% --- Row 1 ---
	\begin{subfigure}[b]{0.4\textwidth}
		\centering
		\includegraphics[width=\textwidth]{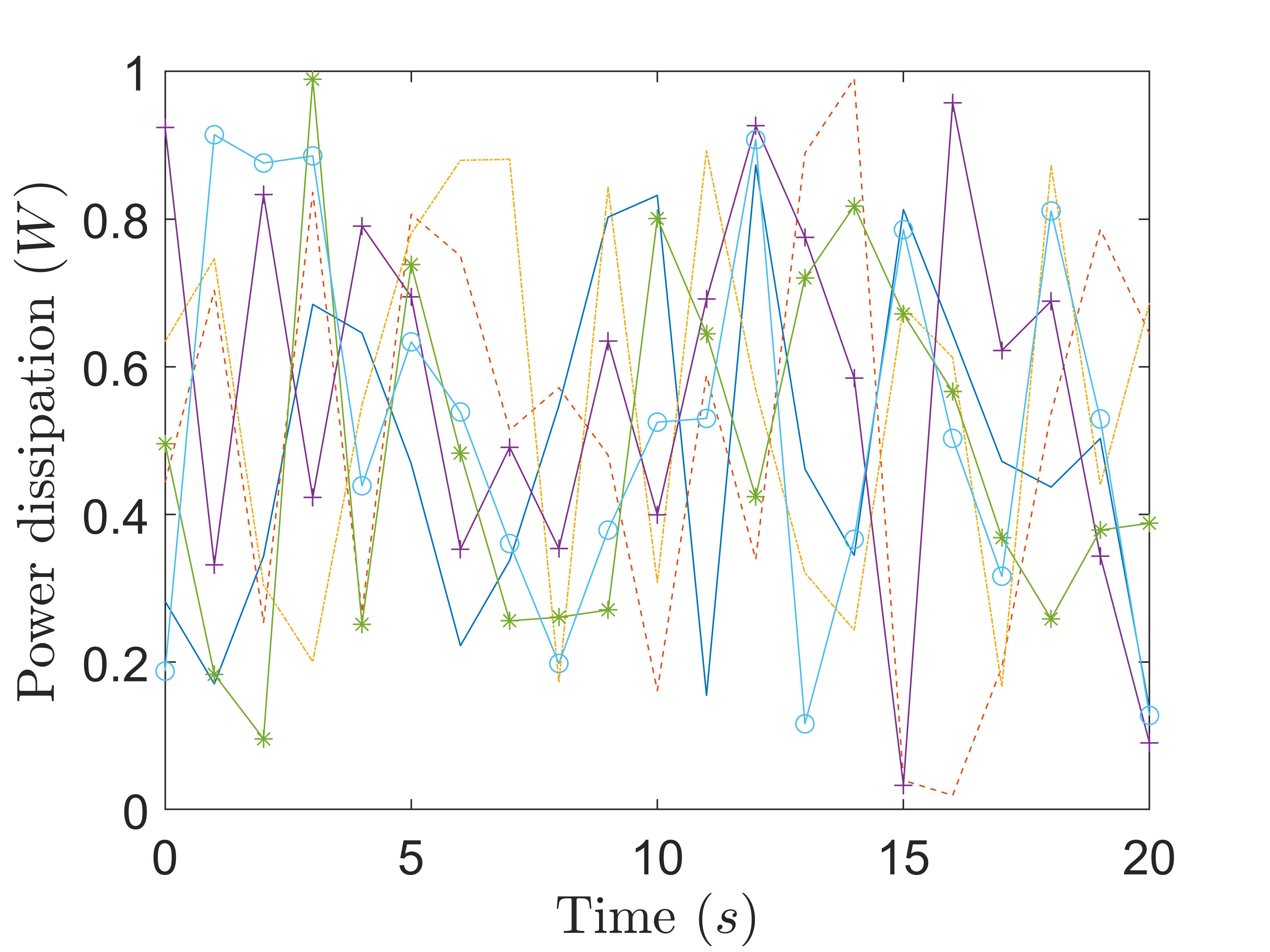}
		\caption{Model development: Input $P_1 (-), P_2 (--), \newline P_3 (-.), P_4 (-+), P_5 (-*), P_6 (-o)$}
		\label{fig:Case3_Training_P}
	\end{subfigure}
	%  \hfill
	\begin{subfigure}[b]{0.4\textwidth}
		\centering
		\includegraphics[width=\textwidth]{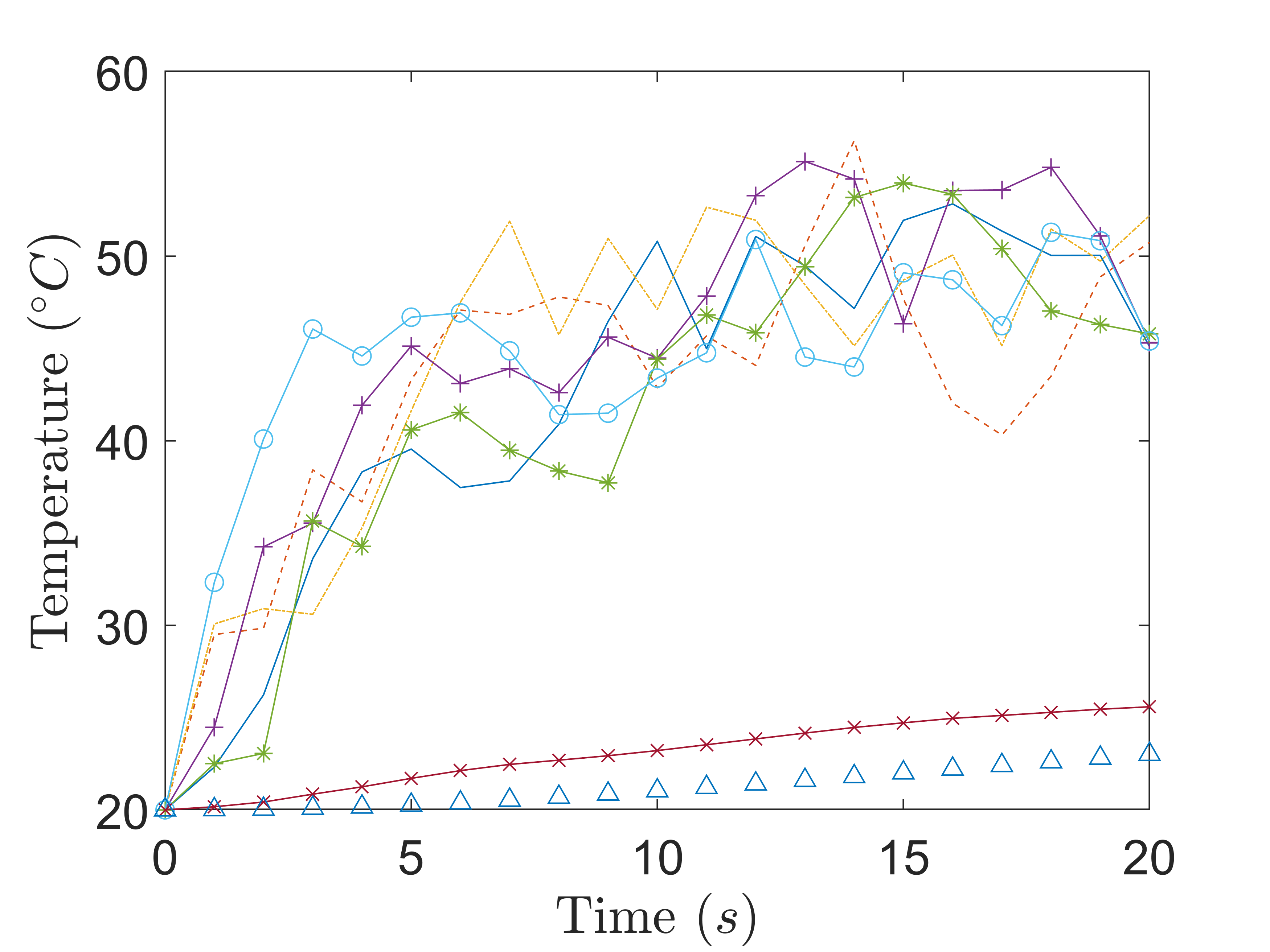}
		\caption{Model development: Output $T_1 (-), T_2 (--), \newline T_3 (-.), T_4 (-+), T_5 (-*), T_6 (-o), T_7 (\texttt{-x}), T_8 (\bigtriangleup)$}
		\label{fig:Case3_Training_T}
	\end{subfigure}
	\par\medskip % safer than \\ inside figure
	\begin{subfigure}[b]{0.3\textwidth}
		\centering
		\includegraphics[width=\textwidth]{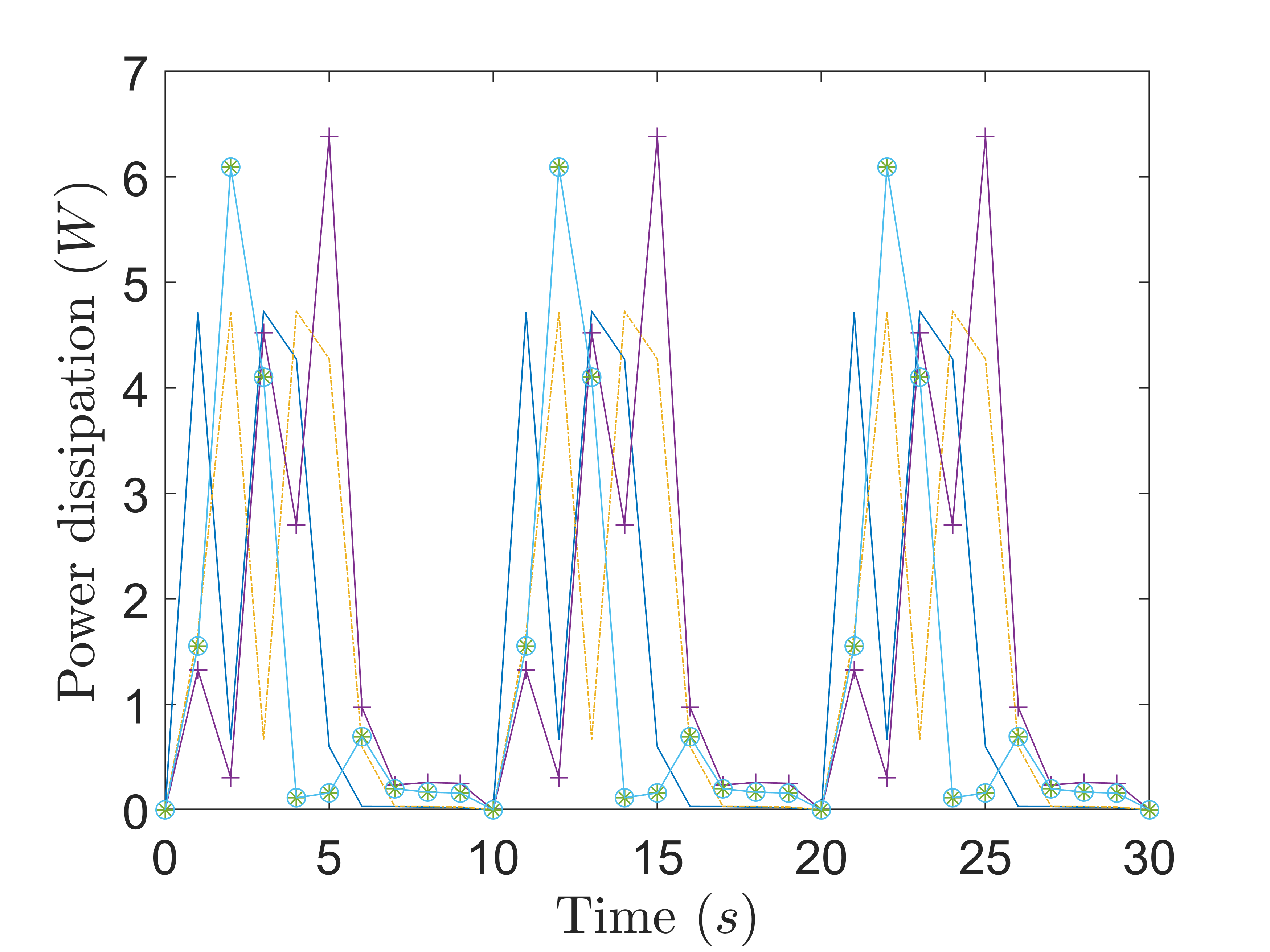}
		\caption{Validation: Input $P_1 (-), P_2 (--), \newline P_3 (-.), P_4 (-+), P_5 (-*), P_6 (-o)$}
		\label{fig:Case3_Test_P}
	\end{subfigure}
	\hfill
	\begin{subfigure}[b]{0.3\textwidth}
		\centering
		\includegraphics[width=\textwidth]{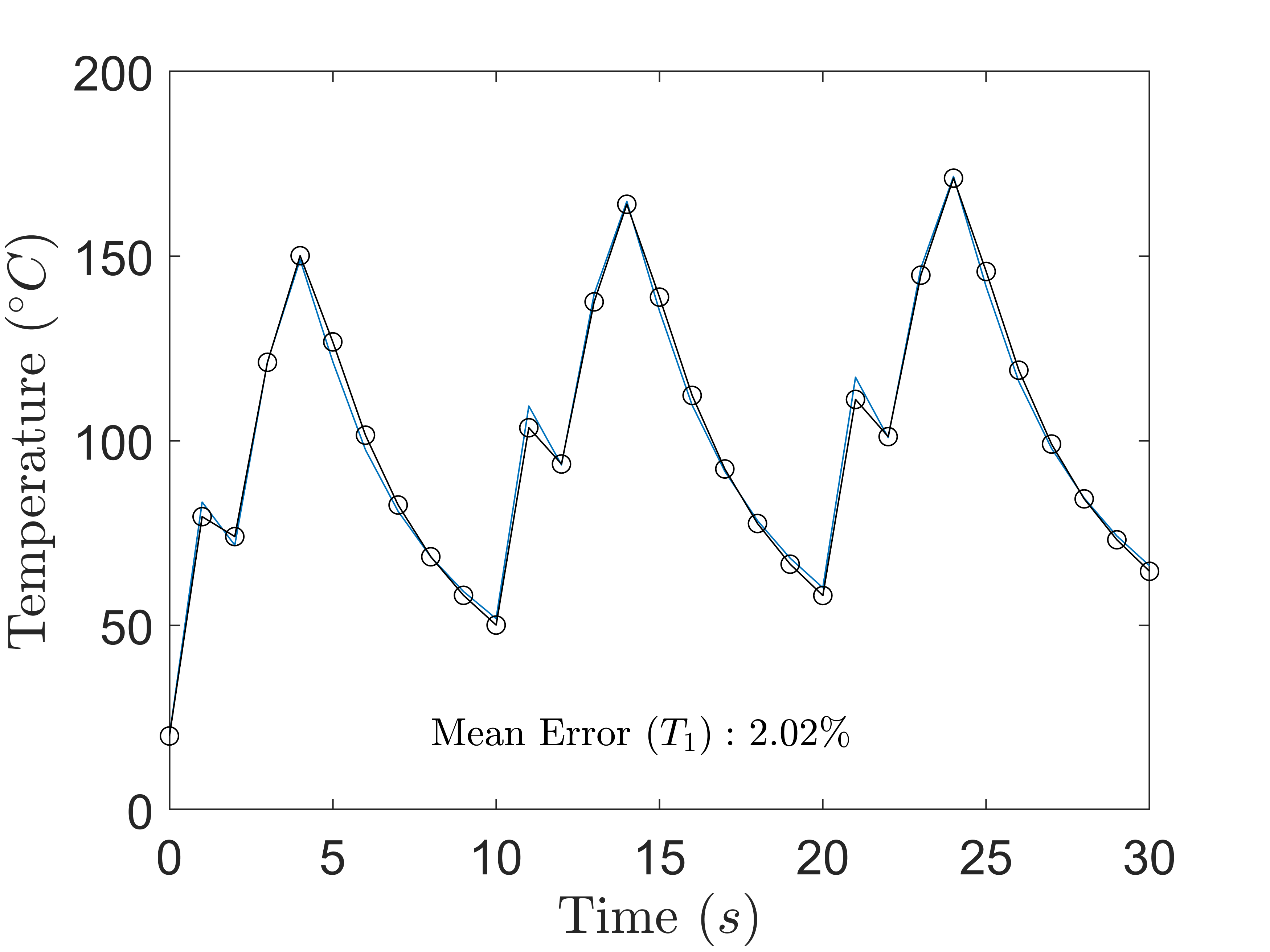}
		\caption{Validation: Simulation $T_1 (-)$ vs model $T_1 (-o)$}
	\end{subfigure}
	\hfill
	\begin{subfigure}[b]{0.3\textwidth}
		\centering
		\includegraphics[width=\textwidth]{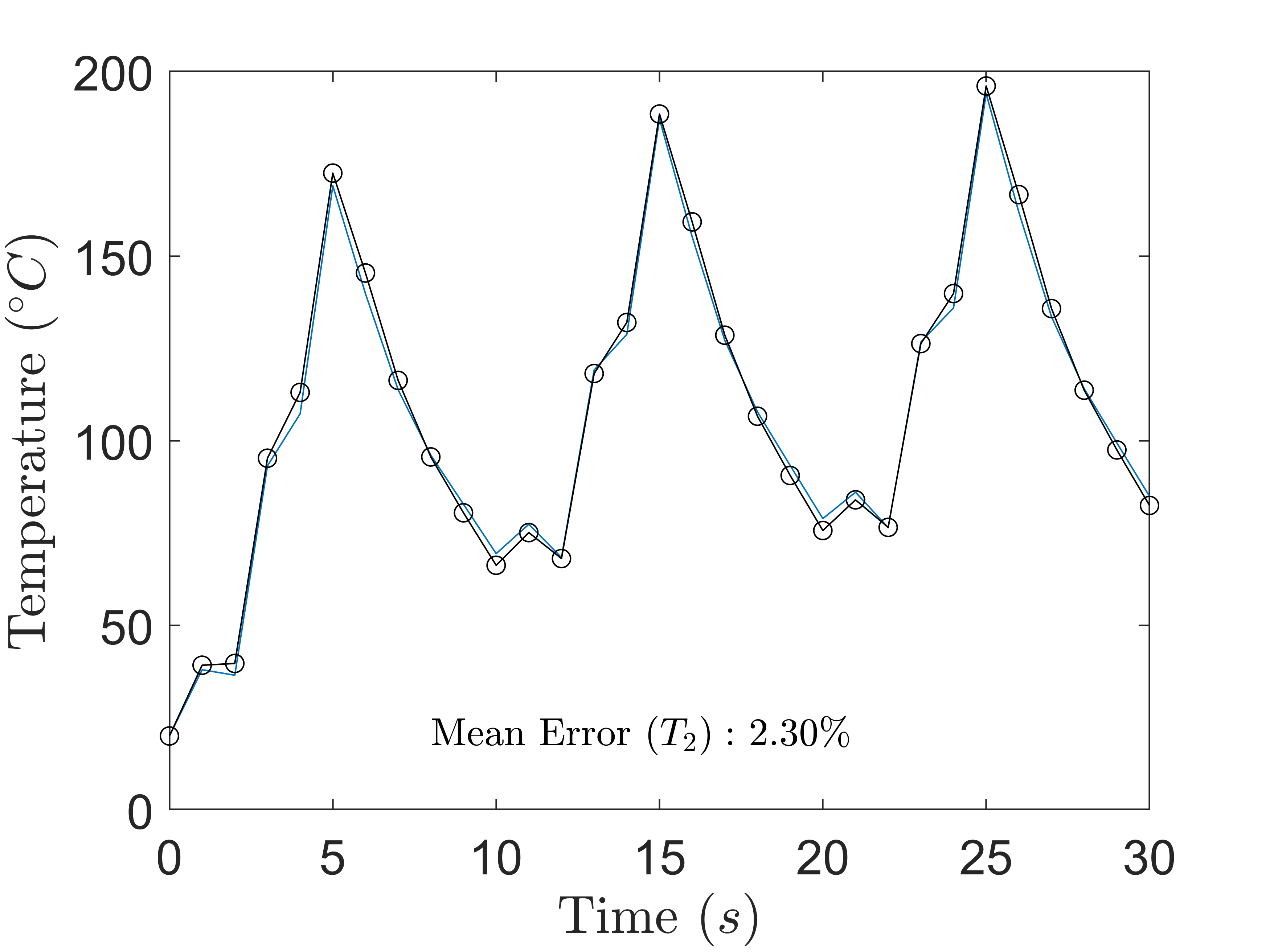}
		\caption{Validation: Simulation $T_2 (-)$ vs model $T_2 (-o)$}
	\end{subfigure}
	\par\medskip % safer than \\ inside figurel
	\begin{subfigure}[b]{0.3\textwidth}
		\centering
		\includegraphics[width=\textwidth]{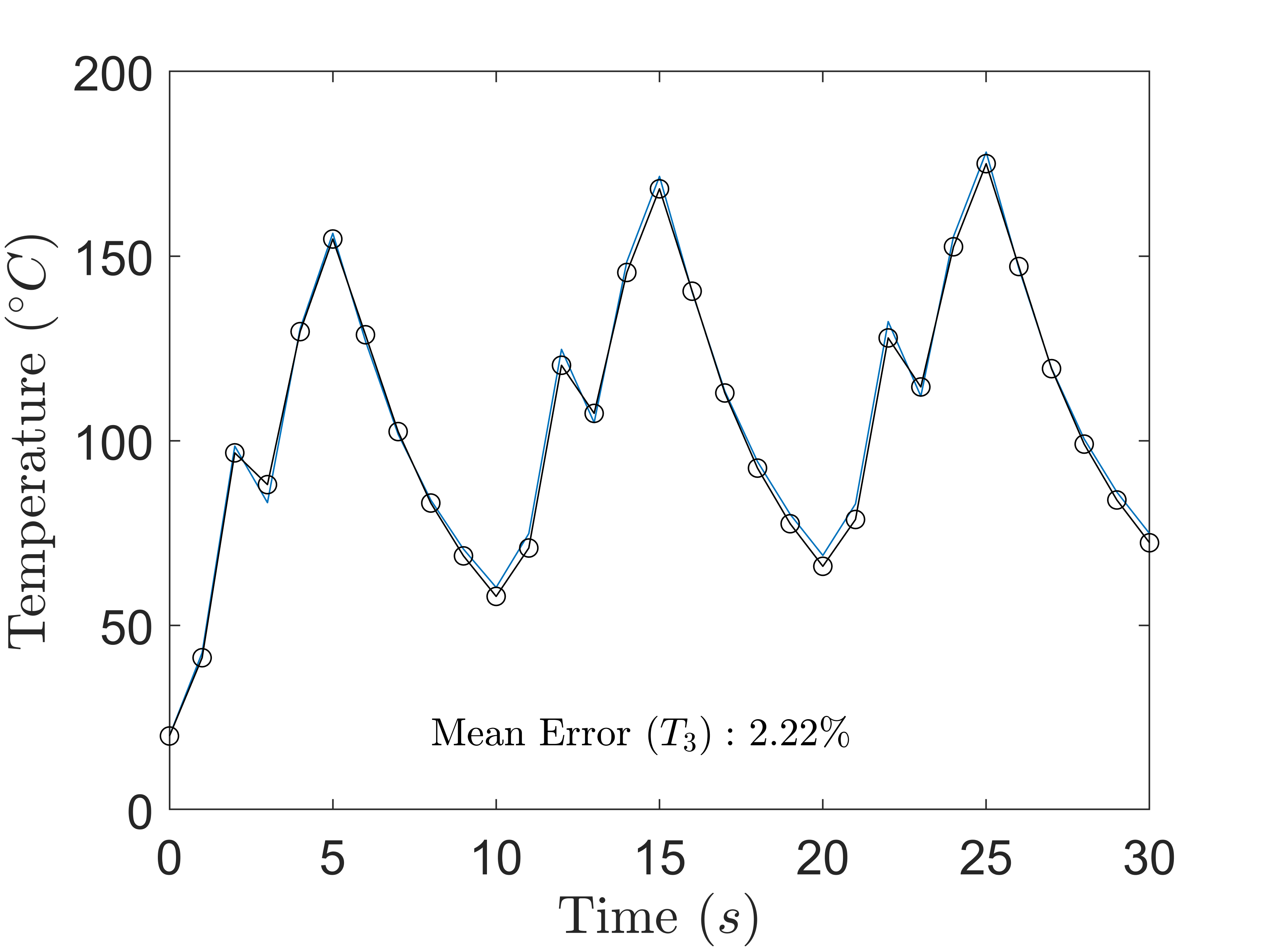}
		\caption{Validation: Simulation $T_3 (-)$ vs model $T_3 (-o)$}
	\end{subfigure}
	\hfill
	\begin{subfigure}[b]{0.3\textwidth}
		\centering
		\includegraphics[width=\textwidth]{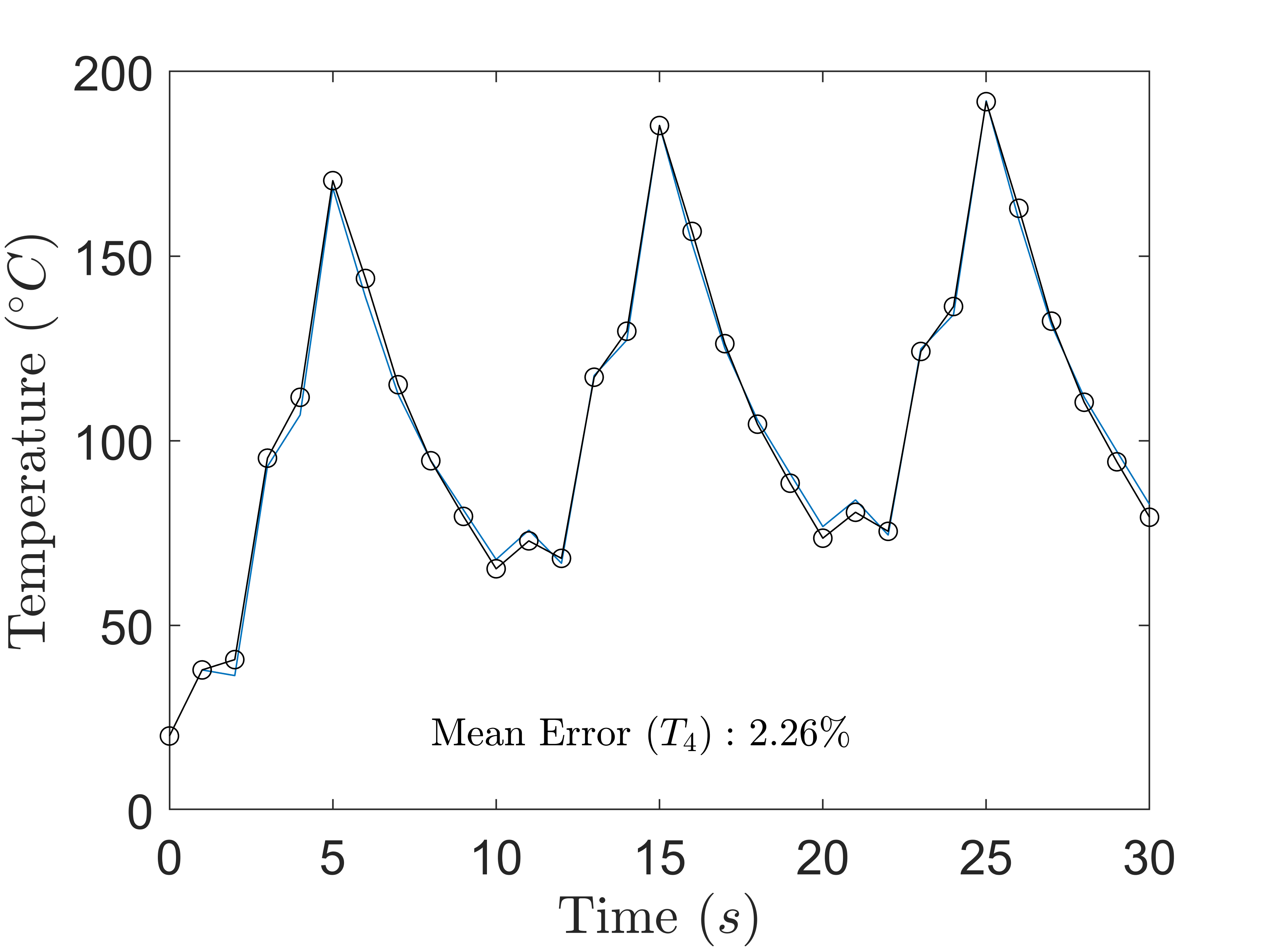}
		\caption{Validation: Simulation $T_4 (-)$ vs model $T_4 (-o)$}
	\end{subfigure}
	\hfill
	\begin{subfigure}[b]{0.3\textwidth}
		\centering
		\includegraphics[width=\textwidth]{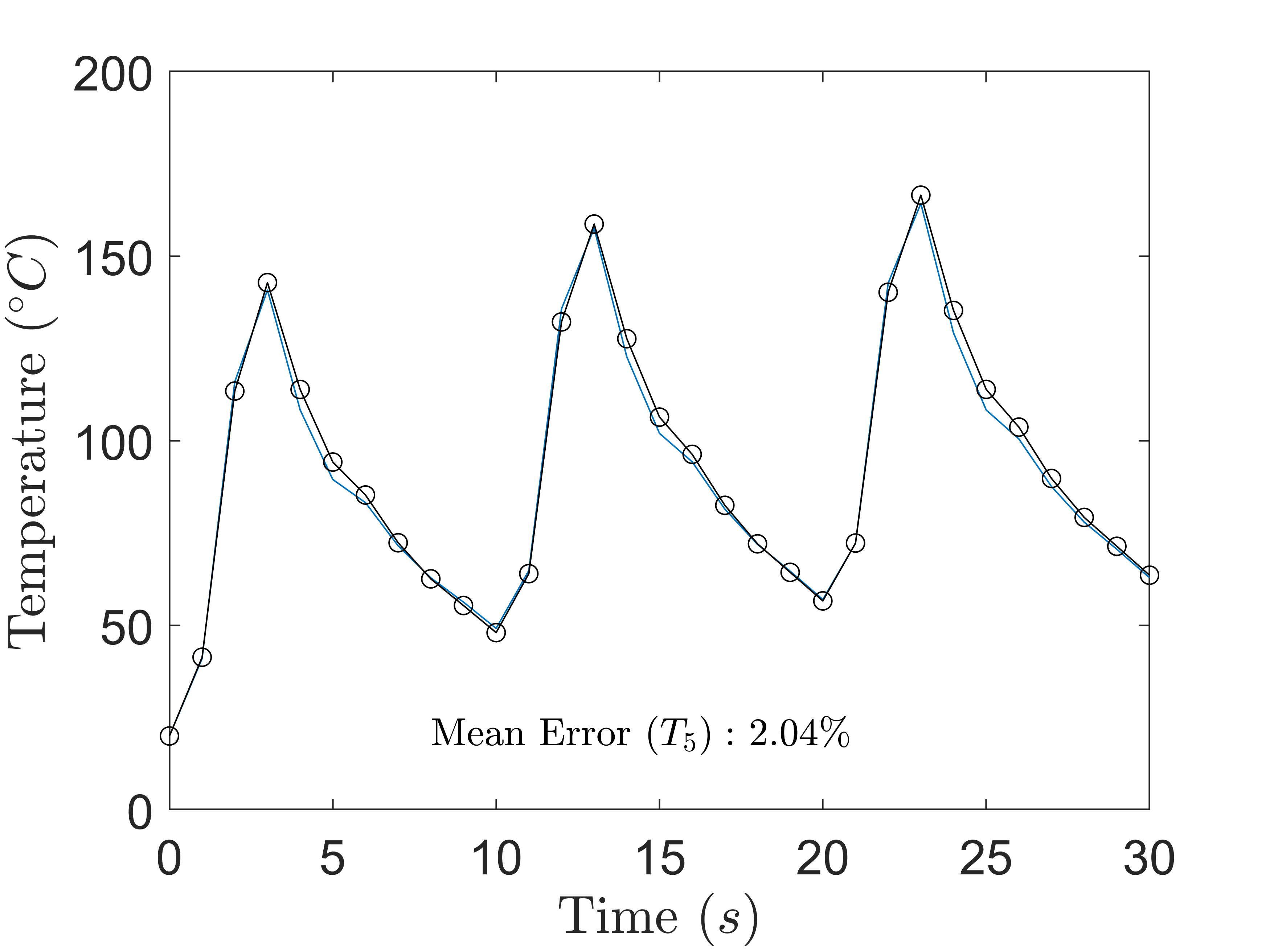}
		\caption{Validation: Simulation $T_5 (-)$ vs model $T_5 (-o)$}
	\end{subfigure}
	\par\medskip % safer than \\ inside figure
	% row 5
	\begin{subfigure}[b]{0.3\textwidth}
		\centering
		\includegraphics[width=\textwidth]{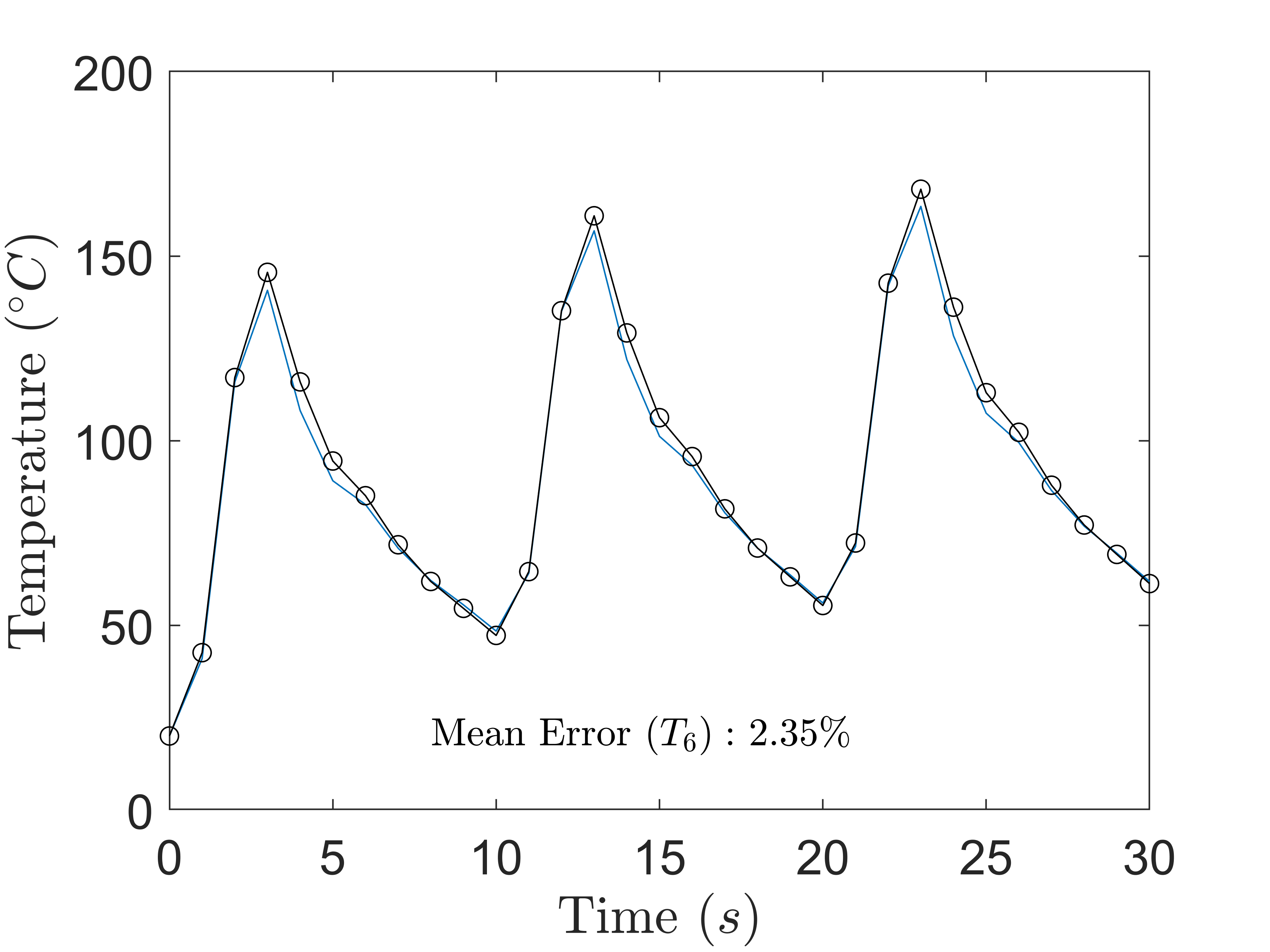}
		\caption{Validation: Simulation $T_6 (-)$ vs model $T_6 (-o)$}
	\end{subfigure}
	\hfill
	\begin{subfigure}[b]{0.3\textwidth}
		\centering
		\includegraphics[width=\textwidth]{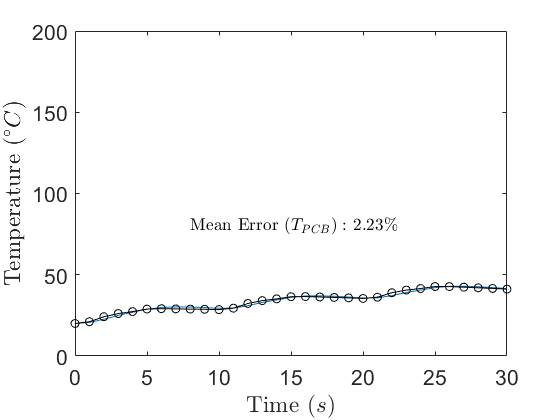}
		\caption{Validation: Simulation $T_7 (-)$ vs model $T_7 (-o)$}
	\end{subfigure}
	\hfill
	\begin{subfigure}[b]{0.3\textwidth}
		\centering
		\includegraphics[width=\textwidth]{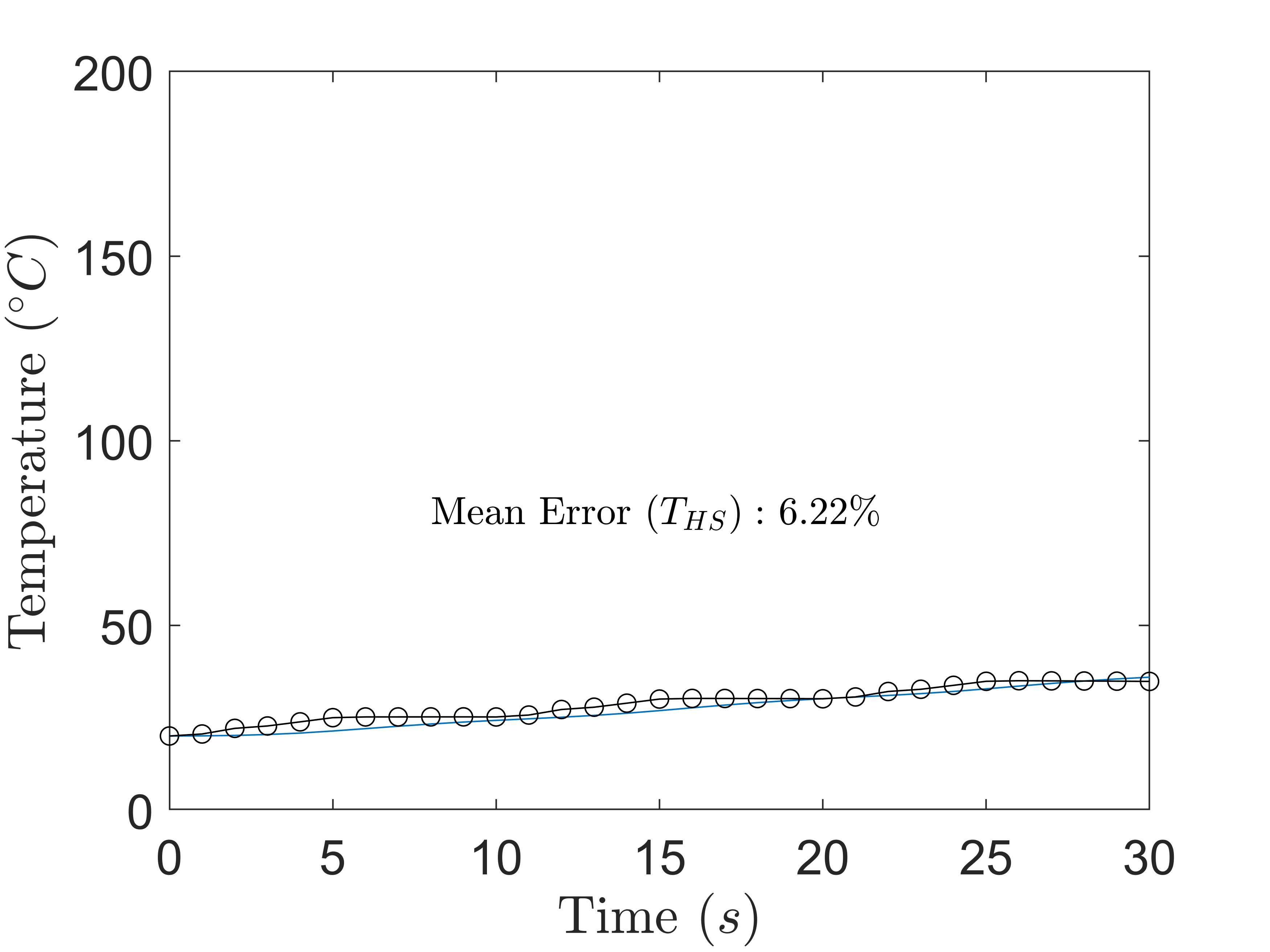}
		\caption{Validation: Simulation $T_8 (-)$ vs model $T_8 (-o)$}
	\end{subfigure}
	\par\medskip % safer than \\ inside figure
	\caption{Comparison of temperatures from simulation and LPLSP model for Inverter module with 6 MOSFETs on a PCB attached to a heatsink, in a natural convection environment. The generation of training / model-development data via CFD requires $600 s$. A traditional parametric study would incur this cost multiplied by the number of heat sources (6 in this case). The mean percentage error between the model predictions and the simulation results is also reported.}
	\label{fig:Case3}
\end{figure}

Comparisons between the temperatures obtained from CFD simulations and those predicted by the LPLSP model using the rank-reduction approach are presented for the two-body conduction case in Fig.\ref{fig:Case1}, the three-body convection case in Fig.\ref{fig:Case2}, and the inverter module in Fig.\ref{fig:Case3}. Accuracy of the predicted temperatures is highly comparable across all algorithms, therefore, for brevity, only the results obtained using the rank-reduction method are presented in the figures. In addition, Table \ref{tab:computation_time} summarizes the computation times associated with model construction and the execution of each model for new transient inputs across the different algorithms proposed in this work. From these comparisons, it is evident that the accuracy of the LPLSP model regardless of the ensemble parameter estimation method employed is consistently high, with all approaches achieving temperature predictions within approximately $5\%$ of the CFD results. Thus, the principal distinction among the methods lies in the computational cost required for model development and subsequent prediction. \\
Algorithm \ref{pc1} is straightforward to construct and performs well for systems with a small number of sources and monitor points (typically $<3$). However, as the system size increases (e.g., more than six coupled nodes), the model development time grows substantially, reaching approximately $t \approx 700\,\text{s}$, which is comparable to the time required for the full CFD simulation itself. It is important to note that this cost is incurred only once during model construction; once the parameters are estimated, the execution time for predicting new transient inputs remains very low, on the order of $O(10^0 s)$. In contrast, Algorithms \ref{pc3} and \ref{pc4} (the rank-reduction and the two-stage estimation with vectorization of core computations) yield significant reductions in model-development time, with values of $t \approx 40 s$ for the rank-reduction method and $t \approx 8 s$ for the two-stage approach for the case of inverter module. This corresponds to an order-of-magnitude reduction in computation time of approximately $O(10^1)$--$O(10^2)$ compared with Algorithm \ref{pc1}. The structure of two-stage method used for the inverter module is similar to that of Algorithm \ref{pc2}, but splitting the estimation in stages dramatically reduces the number of parameters to be estimated. This represents a significant time reduction in both Stages B and C of Fig. \ref{fig:ROMstages}. \\
Finally, the entire workflow from importing the training CFD dataset, performing ensemble parameter estimation, storing the identified parameters, and generating temperature predictions for arbitrary new transient inputs can be fully automated. Once the initial CFD simulations or generation of data from a physical test setup for the training dataset are complete, the remaining steps can be executed rapidly to produce highly accurate reduced-order models.
%----------------
\begin{sidewaystable}[htbp]
	\centering
	\resizebox{\textwidth}{!}{%
		\begin{tabular}{llccc}
			\toprule
			\multirow{2}{*}{\textbf{Case study}} & \multirow{2}{*}{\textbf{Method}} & \multicolumn{2}{c}{\textbf{One time model development cost (s)}} & \textbf{Computation time (s)} \\
			\cmidrule(lr){3-4} \cmidrule(lr){5-5}
			& & \textbf{Simulation for training data} & \textbf{Model development time} & \textbf{Execution time for new inputs} \\
			\midrule
			\multirow{5}{*}{Two-body conduction} 
			& CFD Simulation & - & - & 556 \\
			& Traditional LPLSP (individual parameter estimation) & 324 & 35 & 0.19 \\
			& Algorithm \ref{pc1}: Naive implementation & 162 & 8.283 & 0.19 \\
			& Algorithm \ref{pc2}: Symmetric parameter & 162 & 0.797 & 0.19 \\
			& Algorithm \ref{pc3}: Rank Reduction & 162 & 1.46 & 0.19 \\
			\midrule
			\multirow{5}{*}{Three-body convection} 
			& CFD Simulation & - & - & 610 \\
			& Traditional LPLSP (individual parameter estimation) & 1275 & 90 & 1.12 \\
			& Algorithm \ref{pc1}: Naive implementation & 425 & 39.33 & 1.12 \\
			& Algorithm \ref{pc2}: Symmetric parameter & 425 & 1.334 & 1.12 \\
			& Algorithm \ref{pc3}: Rank Reduction & 425 & 1.90 & 1.12 \\
			\midrule
			\multirow{5}{*}{\makecell{Inverter module \\ (6 MOSFETs, PCBA, Heatsink) \\ - Natural Convection}} 
			& CFD Simulation & - & - & 960 \\
			& Traditional LPLSP (individual parameter estimation) & 3600 & 250 & 2.42 \\
			& Algorithm \ref{pc1}: Naive implementation & 600 & 738.15 & 2.42 \\
			& Algorithm \ref{pc3}: Rank Reduction & 600 & 43.52 & 2.42 \\
			& Algorithm \ref{pc4}: Two-stage parameter estimation & 600 & 8.38 & 2.42 \\
			\midrule
			\multirow{5}{*}{\makecell{Inverter module \\ (6 MOSFETs, PCBA, Heatsink) \\ - Forced Convection $U=10 m/s$}} 
			& CFD Simulation & - & - & 927 \\
			& Traditional LPLSP (individual parameter estimation) & 3750 & 250 & 2.42 \\
			& Algorithm \ref{pc1}: Naive implementation & 625 & 420.2 & 2.42 \\
			& Algorithm \ref{pc3}: Rank Reduction & 625 & 23 & 2.42 \\
			& Algorithm \ref{pc4}: Two-stage parameter estimation & 625 & 9.6 & 2.42 \\
			\bottomrule
		\end{tabular}%
	}
	\caption{Comparison of computation times for CFD simulation and LPLSP model using ensemble parameter estimation with different algorithms. The model development costs for LPLSP model (all algorithms) include the time to run the full CFD simulation to generate the training data and subsequent time to post process the data and create the ROM. For traditional LPLSP method, the model development time is significantly larger because of the number of simulations and post-processing time required for each. This model development time is a one time cost and subsequent execution times for new inputs are of the order $t \approx O(10^0) s$. All executions were performed on a workstation equipped with an Intel\textsuperscript{\textregistered} Core\textsuperscript{\texttrademark} i7-10850H processor (12 cores, 2.7 GHz) and 32 GB RAM.}
	\label{tab:computation_time}
\end{sidewaystable}
%---------------------------2 BODY CONDUCTION 
%--------------------------- 3 BODY CONVECTION
%--------------------------- INVERTER
%-------- INVERTER - FORCED CONVECTION
\begin{figure}[htbp]
	\centering
	% --- Row 1 ---
	\begin{subfigure}[b]{0.4\textwidth}
		\centering
		\includegraphics[width=\textwidth]{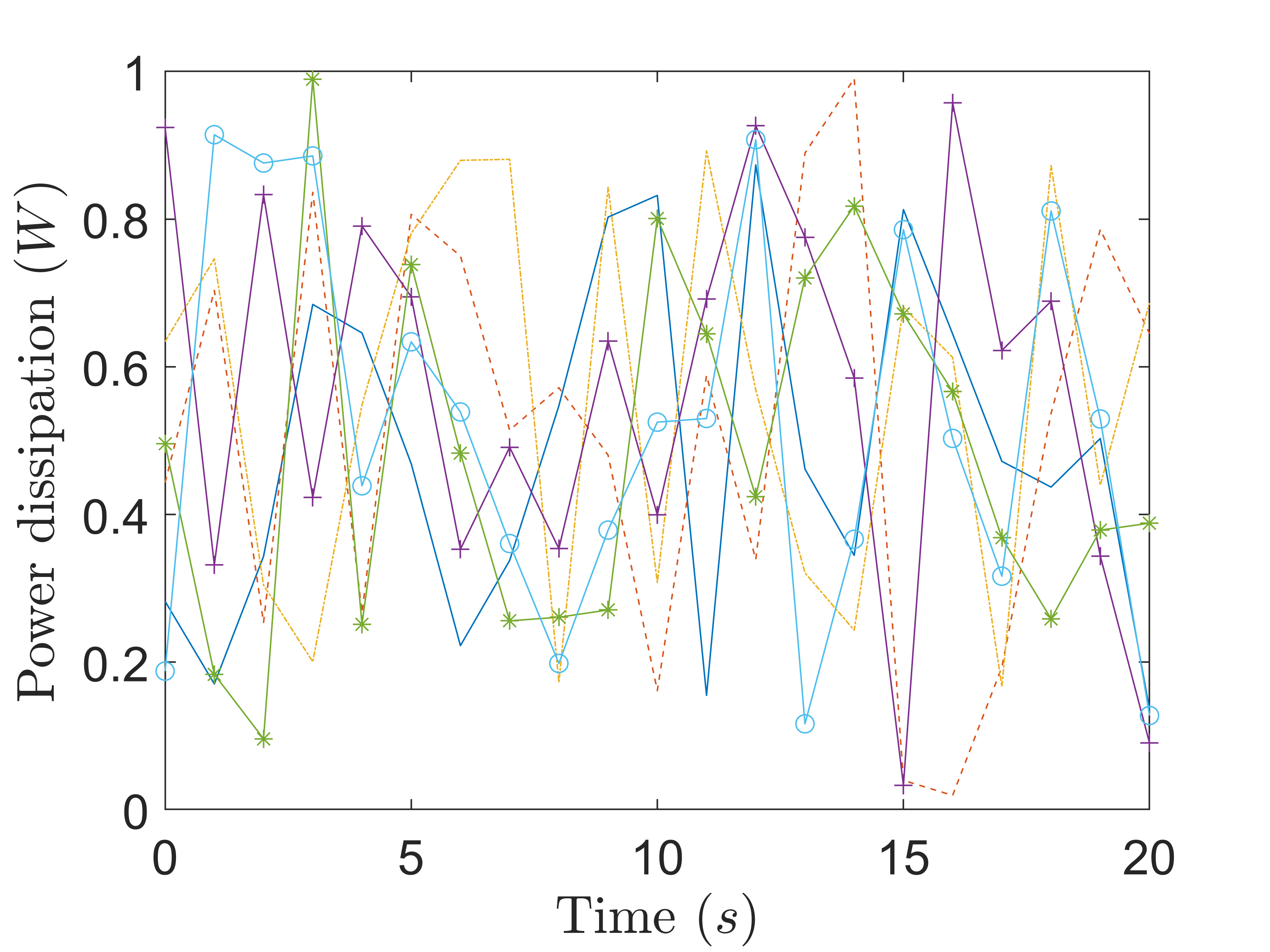}
		\caption{Model development: Input $P_1 (-), P_2 (--), \newline P_3 (-.), P_4 (-+), P_5 (-*), P_6 (-o)$}
		\label{fig:Case4_Training_P}
	\end{subfigure}
	%  \hfill
	\begin{subfigure}[b]{0.4\textwidth}
		\centering
		\includegraphics[width=\textwidth]{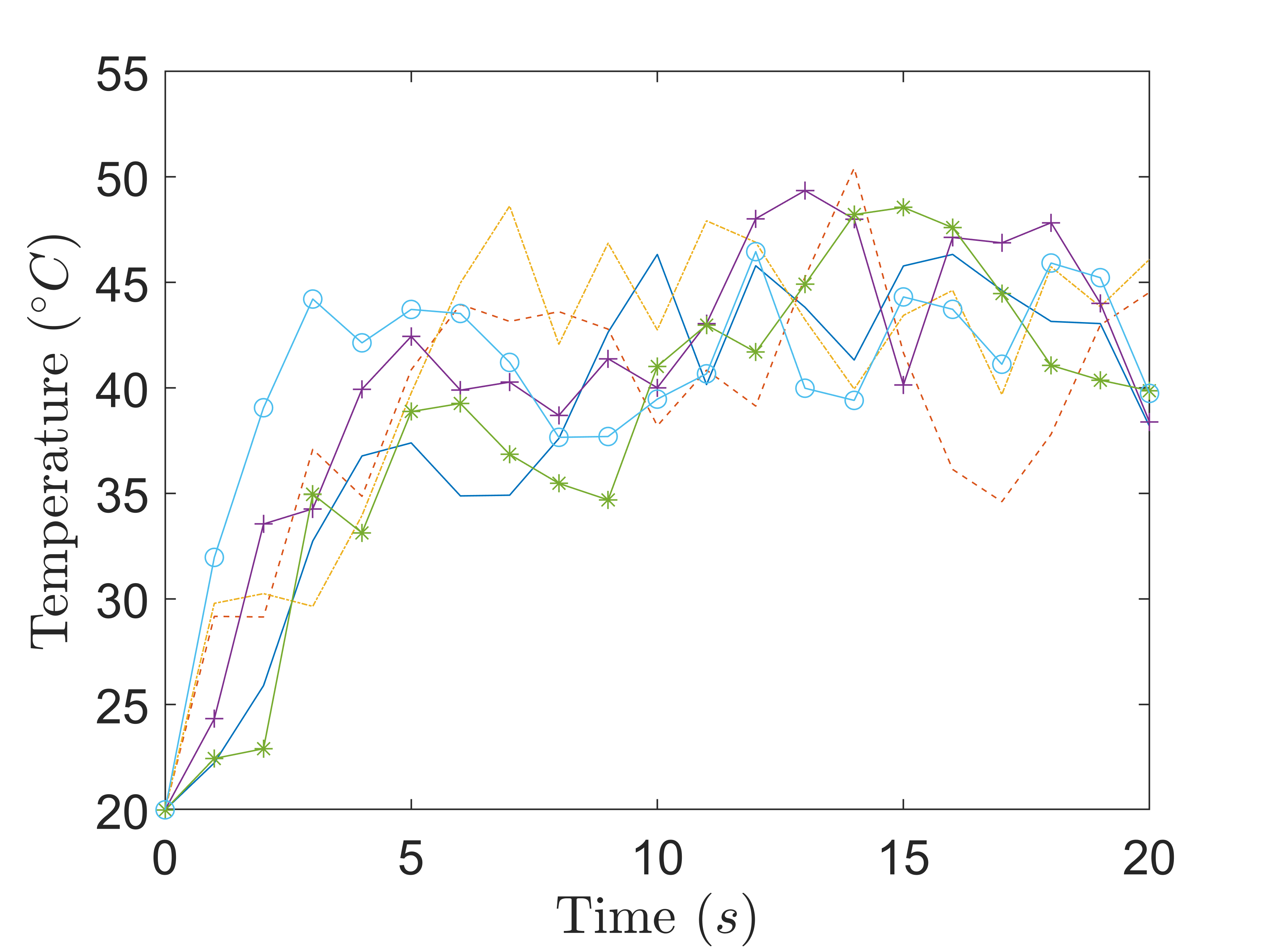}
		\caption{Model development: Output $T_1 (-), T_2 (--), \newline T_3 (-.), T_4 (-+), T_5 (-*), T_6 (-o), T_7 (\texttt{-x}), T_8 (\bigtriangleup)$}
		\label{fig:Case4_Training_T}
	\end{subfigure}
	\par\medskip % safer than \\ inside figure
	\begin{subfigure}[b]{0.3\textwidth}
		\centering
		\includegraphics[width=\textwidth]{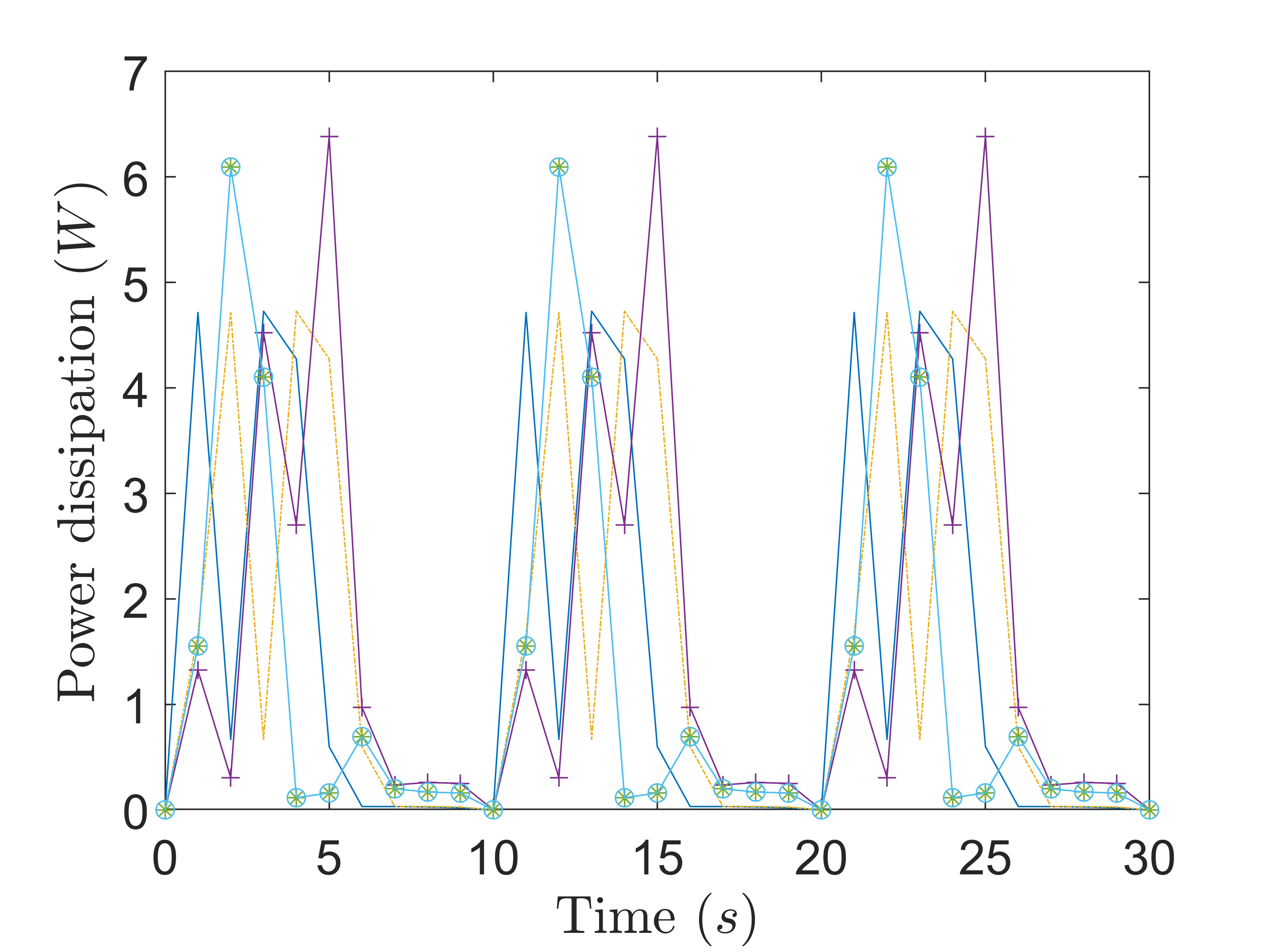}
		\caption{Validation: Input $P_1 (-), P_2 (--), \newline P_3 (-.), P_4 (-+), P_5 (-*), P_6 (-o)$}
		\label{fig:Case4_Test_P}
	\end{subfigure}
	\hfill
	\begin{subfigure}[b]{0.3\textwidth}
		\centering
		\includegraphics[width=\textwidth]{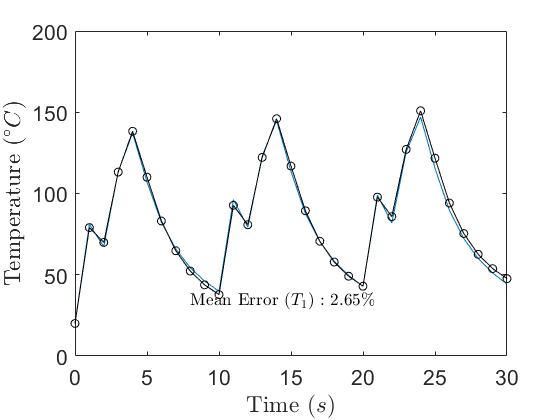}
		\caption{Validation: Simulation $T_1 (-)$ vs model $T_1 (-o)$}
	\end{subfigure}
	\hfill
	\begin{subfigure}[b]{0.3\textwidth}
		\centering
		\includegraphics[width=\textwidth]{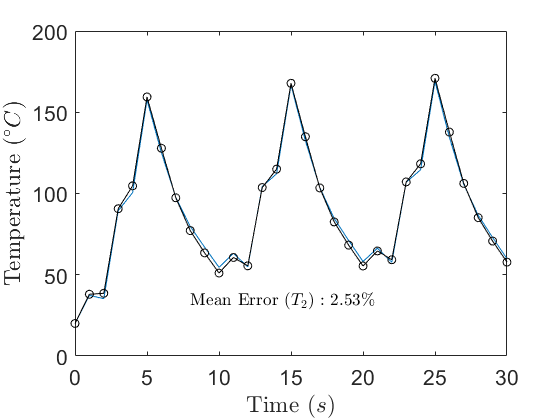}
		\caption{Validation: Simulation $T_2 (-)$ vs model $T_2 (-o)$}
	\end{subfigure}
	\par\medskip % safer than \\ inside figurel
	\begin{subfigure}[b]{0.3\textwidth}
		\centering
		\includegraphics[width=\textwidth]{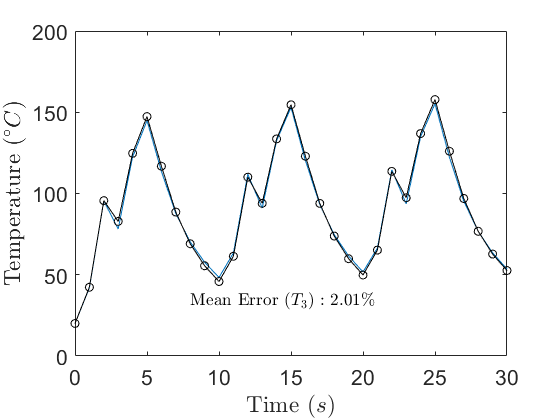}
		\caption{Validation: Simulation $T_3 (-)$ vs model $T_3 (-o)$}
	\end{subfigure}
	\hfill
	\begin{subfigure}[b]{0.3\textwidth}
		\centering
		\includegraphics[width=\textwidth]{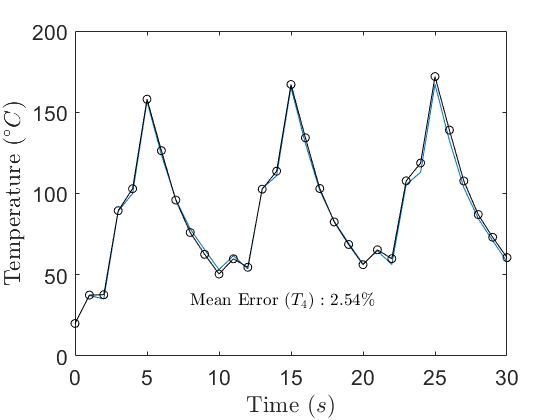}
		\caption{Validation: Simulation $T_4 (-)$ vs model $T_4 (-o)$}
	\end{subfigure}
	\hfill
	\begin{subfigure}[b]{0.3\textwidth}
		\centering
		\includegraphics[width=\textwidth]{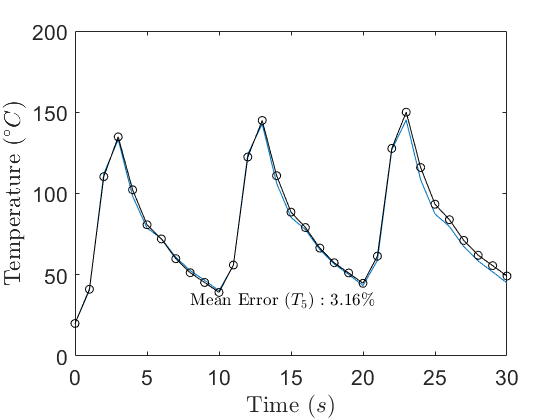}
		\caption{Validation: Simulation $T_5 (-)$ vs model $T_5 (-o)$}
	\end{subfigure}
	\par\medskip % safer than \\ inside figure
	% row 5
	\begin{subfigure}[b]{0.3\textwidth}
		\centering
		\includegraphics[width=\textwidth]{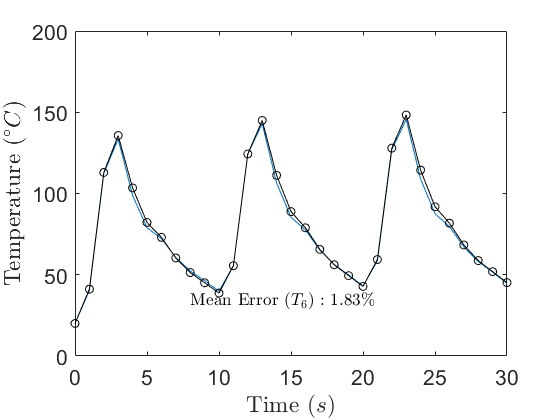}
		\caption{Validation: Simulation $T_6 (-)$ vs model $T_6 (-o)$}
	\end{subfigure}
	\hfill
	\begin{subfigure}[b]{0.3\textwidth}
		\centering
		\includegraphics[width=\textwidth]{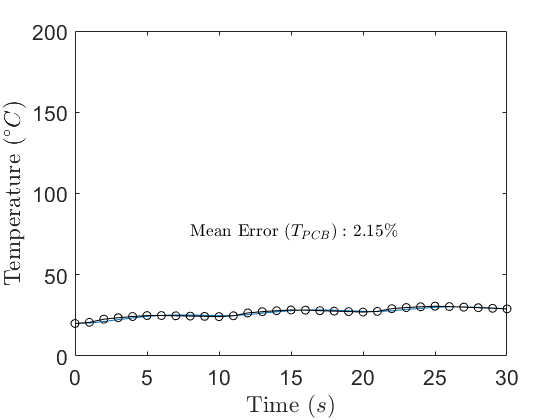}
		\caption{Validation: Simulation $T_7 (-)$ vs model $T_7 (-o)$}
	\end{subfigure}
	\hfill
	\begin{subfigure}[b]{0.3\textwidth}
		\centering
		\includegraphics[width=\textwidth]{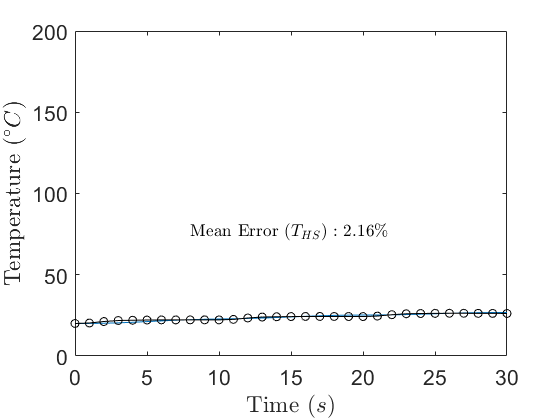}
		\caption{Validation: Simulation $T_8 (-)$ vs model $T_8 (-o)$}
	\end{subfigure}
	\par\medskip % safer than \\ inside figure
	\caption{Comparison of temperatures from simulation and LPLSP model for Inverter module with 6 MOSFETs on a PCB attached to a heatsink, in a forced convection environment with $U=10 m/s$. The generation of training / model-development data via CFD requires $625 s$. A traditional parametric study would incur this cost multiplied by the number of heat sources (6 in this case). The mean percentage error between the model predictions and the simulation results is also reported.}
	\label{fig:Case4}
\end{figure}
%%%%%%%%%%%%%%%%%%%%%%%
\section{Conclusion}
This work presents an improved workflow for constructing reduced-order thermal models using the Lumped Parameter Linear Superposition (LPLSP) framework. By replacing the traditional sequence of parametric studies with an ensemble parameter estimation approach, the full set of model coefficients can be obtained from a single transient dataset. This fundamentally streamlines the model-development process and removes practical barriers associated with isolating individual heat sources, particularly in hardware environments. The proposed parameter estimation strategies including rank-reduced method and two-stage approach produce models that retain the accuracy characteristic of the original LPLSP method while significantly reducing computation time. The resulting ROMs enable rapid evaluation of new transient operating conditions and can be generated with a high degree of automation, making the approach well suited for CFD-based model development and for constructing thermal digital twins of physical systems. This has tremendous applications in a multitude of industries including automotive, aerospace and healthcare. This study primarily presents results for conduction, natural and forced convection conditions. Future work will extend the method to systems with forced convection with variable flow velocities, as well as to the creation of digital twins for physical test setups involving multiple modes of heat transfer including radiation.

\FloatBarrier

\printbibliography
%\bibliographystyle{unsrt}
%\bibliography{references}
\end{document}